\documentclass{article}
\usepackage[letterpaper,top=2cm,bottom=2cm,left=2cm,right=2cm,marginparwidth=1.75cm]{geometry}

\usepackage{graphicx} 
\usepackage{xcolor}
\usepackage{enumitem}
\usepackage{orcidlink}
\usepackage{textgreek}
\usepackage{calrsfs}
\usepackage{amssymb}
\usepackage{amsmath}
\usepackage{amsthm}
\usepackage{amsfonts}
\usepackage{subcaption}
\usepackage{multirow}
\usepackage{graphicx}
\usepackage{amsfonts}
\usepackage{mathtools}
\usepackage{mathrsfs}
\usepackage{placeins}
\usepackage{multirow}
\usepackage{verbatim,hyperref}
\usepackage{xcolor}
\usepackage{stmaryrd}
\usepackage{array}
\usepackage{caption}
\newcolumntype{M}[1]{>{\centering\arraybackslash}m{#1}}
\usepackage{algorithm}
\usepackage{algpseudocode}
\usepackage{dirtytalk}
\usepackage{todonotes}
\usepackage{authblk}
\usepackage{svg}
\usepackage{subcaption}
\usepackage{listings}
\usepackage{filecontents}
\usepackage{tikz}
\usepackage{pgfplots}
\usepackage{tabularx}

\usepackage{ulem}

\newcommand{\averagel}{\{\!\!\{}
\newcommand{\averager}{\}\!\!\}}

\newcommand{\jumpl}{[\![}
\newcommand{\jumpr}{]\!]}

\newcommand{\partition}{\mathcal{T}_h}
\newcommand{\facesinternal}{\mathcal{F}^\mathrm{I}_h}

\newcommand{\facesboundary}{\mathcal{F}^\mathrm{B}_h}
\newcommand{\facesboundaryvent}{\mathcal{F}^\mathrm{B,Vent}_h}
\newcommand{\facesboundarypial}{\mathcal{F}^\mathrm{B,Pial}_h}

\newcommand{\Vh}{V_h^\mathrm{dG}}

\newtheorem{theorem}{Theorem}[section]
\newtheorem{remark}[theorem]{Remark}

\title{A whole-brain model of amyloid beta accumulation and cerebral hypoperfusion in Alzheimer's disease}

\author[1,2]{Mattia Corti \orcidlink{0000-0002-7014-972X}}

\author[3]{Andrew Ahern \orcidlink{0009-0007-3039-4564}}

\author[3]{Alain Goriely \orcidlink{0000-0002-6436-8483}}

\author[4,5]{Ellen Kuhl \orcidlink{0000-0002-6283-935X}}

\author[1]{Paola F. Antonietti \orcidlink{0000-0002-2138-3878}}

\affil[1]{MOX-Dipartimento di Matematica, Politecnico di Milano, Milan, Italy}

\affil[2]{Faculty of Mathematics, University of Vienna, Oskar-Morgenstern-Platz 1, 1090 Vienna, Austria}

\affil[3]{Mathematical Institute, University of Oxford, Oxford, UK}

\affil[4]{Institute of Applied Mechanics, Friedrich-Alexander-Universität Erlangen-Nürnberg, Erlangen, Germany}

\affil[5]{Department of Mechanical Engineering, Stanford University, Stanford, CA, USA}

\date{}

\begin{document}

\maketitle
\begin{abstract}
Accumulation of amyloid beta  proteins is a defining feature of Alzheimer's disease, and is usually accompanied by cerebrovascular pathology. Evidence suggests that amyloid beta and cerebrovascular pathology are mutually reinforcing; in particular, amyloid beta suppresses perfusion by constricting capillaries, and hypoperfusion promotes the production of amyloid beta. Here, we propose a whole-brain model coupling amyloid beta and blood vessel through a hybrid model consisting of a  reaction--diffusion system for the protein dynamics and  porous--medium model of blood flow within and between vascular networks: arterial, capillary and venous. We discretize the resulting parabolic--elliptic system of PDEs by means of a high-order discontinuous Galerkin method in space and an implicit Euler scheme in time. Simulations in realistic brain geometries demonstrate the emergence of multistability, implying  that a sufficiently large pathogenic protein seeds is necessary to trigger disease outbreak. Motivated by the `two-hit vascular hypothesis' of Alzheimer's disease that hypoperfusive vascular damage triggers amyloid beta pathology, we also demonstrate that localized hypoperfusion, in response to injury, can destabilize the healthy steady state and trigger brain-wide disease outbreak.
\end{abstract}

\section{Introduction}
\label{sec:intro}
\begin{figure}[t!]
	\centering
	{\includegraphics[width=\textwidth]{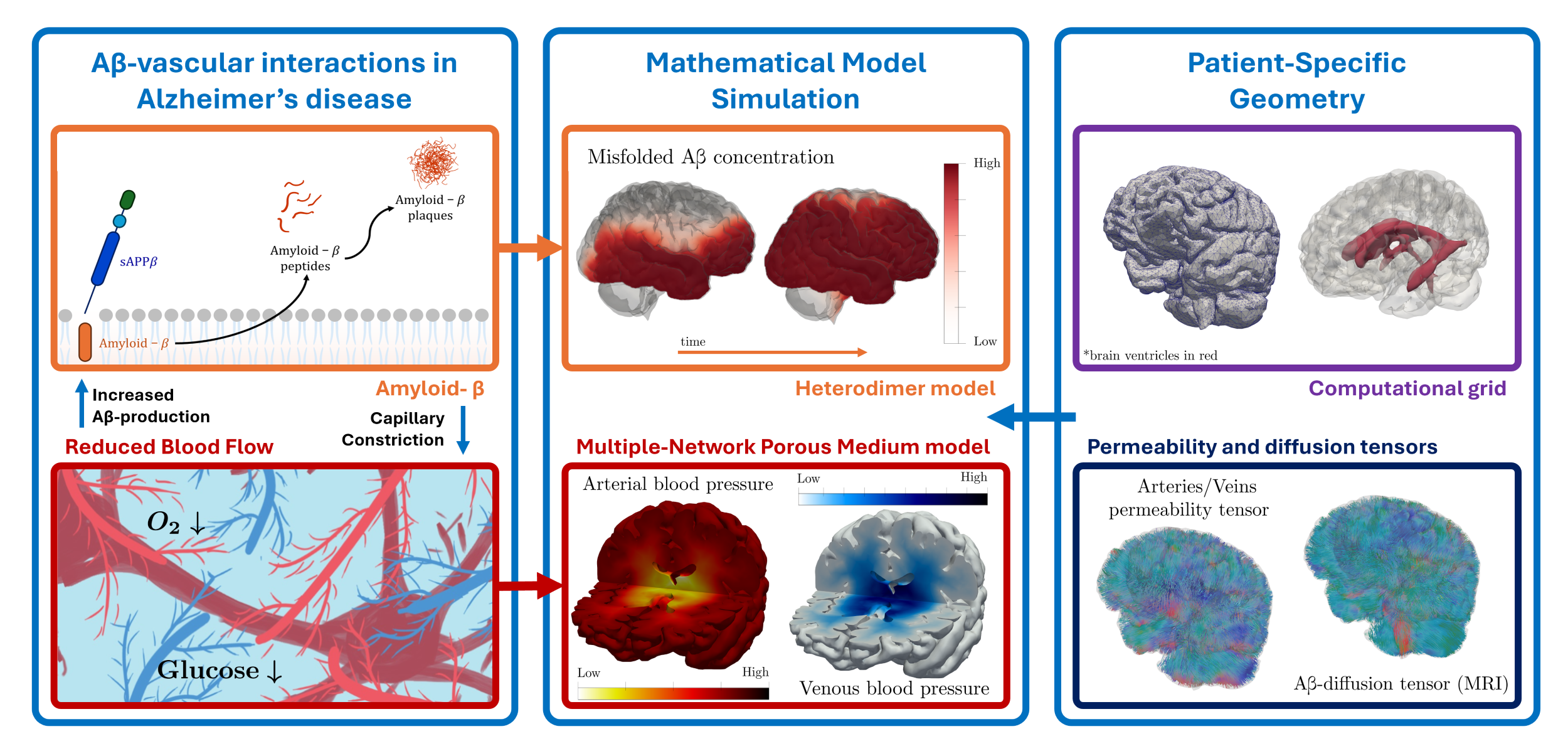}}
	\caption{Synthetic representation of the article structure. Description of the A\textbeta{--}vascular interaction in AD (left panel), patients-specific geometry (right panel), and mathematical models and resulting numerical simulations (center panel).}
	\label{fig:introfig}
\end{figure}
\noindent
\noindent
Alzheimer's disease (AD) is a progressive neurodegenerative disorder characterized by neuronal loss and impaired synaptic communication. It has long been established that the spread of misfolded, prion-like proteins plays a central role in AD and related disorders \cite{scheltens_alzheimers_2021}. In particular, disease onset is associated with the accumulation of two pathological proteins: amyloid-beta (A\textbeta{}) and tau \cite{goedert_alzheimer_2015}. Cerebrovascular abnormalities frequently accompany AD pathology—for instance, cerebral amyloid angiopathy occurs in more than $90\%$ of AD cases \cite{iadecolar_pathobiology_2013}, and reductions in cerebral blood flow (CBF) represent one of the earliest measurable biomarkers of the disease \cite{iturria-medina_early_2016}.

The interplay between A\textbeta{} and cerebral blood flow is well documented. A\textbeta{} acts as a vasoconstrictor, increasing vascular resistance and thereby reducing CBF \cite{thomas_beta_1996,Niwa2001}. Oligomeric A\textbeta{} induces the release of vasoconstrictive agents such as endothelin-1 and generates oxidative stress, which can result in capillary occlusion \cite{cruz_hernandez_neutrophil_2019,nortley_amyloid_2019}. Conversely, hypoxia accelerates amyloid precursor protein (APP) processing and impairs its clearance across the blood–brain barrier, promoting A\textbeta{} accumulation \cite{korte_cerebral_2020,sun_hypoxia_2006}. Hypoperfusion-induced injury can further exacerbate APP aggregation, possibly as a compensatory mechanism \cite{shi_hypoperfusion_2000,hefter_APP_2017}. Collectively, these mechanisms give rise to a positive feedback loop between A\textbeta{} buildup and vascular dysfunction \cite{Iadecola2004,Kalaria2012,korte_cerebral_2020}. Several studies even suggest that vascular damage may act as an early trigger for A\textbeta{} pathology in AD \cite{Zlokovic2011,kisler_cerebral_2017}.

\paragraph{Mathematical models.} In recent years, several mathematical models for the dynamics of \textit{prion-like} proteins have been proposed. The description of the phenomena requires models based on partial differential equations (PDEs) able to describe both the temporal and the spatial dynamics of the phenomena, typically at organ scale. One example developed to describe the dynamics of A\textbeta{} is the Smoluchowski-type models \cite{bertsch_alzheimers_2017,sveva2}, which are  derived from microscopical interactions and distinguishes different sizes of pathological aggregates \cite{franchi_from_2016}. Some simplified models focus on the interactions between healthy and misfolded proteins, such as the heterodimer model \cite{fornari_prion-like_2019,matthaus_diffusion_2006}, or only on the misfolded proteins population, such as the Fisher-Kolmogorov equation \cite{weickenmeier_physics_2019}. These simplified yet informative models offer the key advantage that they can be  fully validated against existing data, while also accurately predicting disease progression at both the personalized and regional scales~\cite{chaggar2025personalised}. Furthermore, the much smaller number of physical parameters relative to Smoluchowski-type models makes the calibration procedure significantly more tractable \cite{corti_uncertainty_2024,corti_exploring_2024}. They have also been extended to describe connections with other physical phenomena occurring in AD. For example, the elastic deformations due to atrophy \cite{schafer_interplay_2019,weickenmeier_multiphysics_2018,pederzoli_coupled_2025}, the clearance mechanisms \cite{brennan_role_2024},  the interactions with tau proteins \cite{thompson_proteinprotein_2020}, and brain activity \cite{goriely2020neuronal,alexandersen2023multi}. 
\par

However, little attention has been devoted to the development of models coupling  A\textbeta{} and CBF. In \cite{ahern_modelling_2025}, the authors propose a modification of the heterodimer model on networks which takes into account CBF reduction. The first goal of the current work is to build on this idea by developing a  continuous PDE model of A\textbeta{} accumulation in the brain describing also the interactions with the cerebral microvasculature through (a) vasoconstriction and (b) perfusion-dependent A\textbeta{} production and clearance (see Figure \ref{fig:introfig}). Our model is based on the heterodimer model for the A\textbeta{}-dynamics description and on a multiple-network porous-medium model for the blood perfusion.
\par

Porous-medium models have been extensively used in the context of brain poromechanics to study strokes \cite{jozsa_porous_2021,josza_mri_2023,jozsa_sensitivity_2021}. Moreover, a poroelastic version of those models have been proposed to study the brain hemodynamics on the heartbeat scale \cite{corti_numerical_2023,lee_mixed_2019,piersanti_parameter_2021,tully_cerebral_2011}, and  clearance mechanisms \cite{li_computational_2023,vardakis_exploring_2020}. 

In our model, we analyze the impact of  A\textbeta{} on the blood flow in the capillaries, altering transfer between  compartments and capillaries permeability. 
\par

\paragraph{Numerical methods and patient-specific simulations.} The potential clinical utility of the constructed mathematical models strongly depend on the design quality of numerical methods to simulate PDEs in patient-specific geometries. In the context of proteinopathies, connectome graph-based simulations for the models introduced in the previous section  have been fully implemented \cite{ahern_modelling_2025,brennan_role_2024,corti_uncertainty_2024,fornari_prion-like_2019,thompson_proteinprotein_2020}. However, the construction of a graph-based model for our purposes would suffer from the limitation of describing the brain’s microvasculature as a collection of small independent capillary networks, losing the coupling with the complete vasculature structure at the organ level \cite{ahern_modelling_2025}.
\par
To study physical effects in the brain, multiple works have introduced numerical discretizations of PDEs in brain geometries based on continuous finite elements \cite{causemann_human_2022,lee_mixed_2019, piersanti_parameter_2021,vardakis_exploring_2020} and discontinuous Galerkin (dG) methods \cite{antonietti_discontinuous_2024,corti_discontinuous_2023,corti_numerical_2023,corti_structure_2024,pederzoli_coupled_2025}. The advantage of using these high-fidelity discretizations is the possibility of constructing patient-specific computational grids from medical images and taking into account geometrical and functional information in the PDE model.
\par
For our mathematical model, we propose a discretization based on a discontinuous Galerkin method in space and implicit Euler time stepping in time. The dG methods provide numerous advantages in our context because they are designed for high-order approximations that are of primary importance for the heterodimer system, which typically admits travelling-wave solutions \cite{kevrekidis2020anisotropic,antonietti_discontinuous_2024,corti_exploring_2024}. Indeed, high-order dG schemes exhibit favourable properties for wave-like problems, allowing one to capture travelling fronts with reduced numerical dispersion and dissipation errors compared with standard low-order continuous finite element approximations as discussed in \cite{antonietti_high_2016,antonietti_high-order_2018}. Moreover, as a possible extension of the present discretization, the dG framework naturally allows for a local adjustment of the polynomial degree $\ell$, which can be exploited to maintain a low computational cost through $\ell$-adaptivity \cite{leimer_p-adaptive_2025}. Additionally, they allow the use of polygonal mesh elements which is particularly useful to reduce computational costs in brain applications. Indeed, using mesh agglomeration strategies \cite{antonietti_polytopal_2026}, the complex boundary and internal interfaces can be accurately described with a low number of mesh elements. Concerning the time discretization, the implicit Euler method is combined with a loosely-coupled strategy to solve the parabolic (heterodimer model) and elliptic (multiple-network porous-medium model) parts of the system separately.
\par
In this work, we first carry out simulations in simplified geometries to study and understand the fundamental properties of our mathematical model. Then we perform realistic numerical simulations in brain geometries. In particular, we show the existence of multiple stable states starting from different levels of initial seeding of A\textbeta{}. Finally, we evaluate the injury-induced initiation of AD pathology, by imposing different levels of initial hypoperfusion in the frontal lobe.
\paragraph{Structure of the manuscript.}  Section~\ref{sec:math_model} is dedicated to the mathematical development of our model. We start from the presentation of the heterodimer model for A\textbeta{} spreading in Section~\ref{sec:HM:model} and the multiple-network porous-medium model in Section~\ref{sec:model:MP}. Then, in Section~\ref{sec:coupled:model}, we propose a novel coupled mathematical model to describe the connections between A\textbeta{} and CBF. In Section~\ref{sec:discretization} we introduce the discretization of the problem based on a dG-space discretization and implicit Euler time stepping. In Section~\ref{sec:num_val}, we present numerical simulations in simple geometries to validate the numerical solver and analyze properties of the mathematical model. In Section~\ref{sec:num_brain}, we simulate the system in realistic three-dimensional brain geometries to study both the multistability of the system depending on the magnitude of the initial seeding and the hypoperfusion-induced pathology. Finally, in Section~\ref{sec:conclusion}, we conclude and discuss further developments.
\section{Mathematical model}
\label{sec:math_model}

\subsection{Heterodimer model}
\label{sec:HM:model}
We adapt Prusiner's heterodimer prion model \cite{prusiner_molecular_1991} for the protein kinetics of A\textbeta{}. That is, we assume that A\textbeta{} monomers can be either normal or pathogenic, with concentrations $u=u(\boldsymbol{x},t)$ and $\tilde{u}=\tilde{u}(\boldsymbol{x},t)$, respectively. Pathogenic monomers can bind to normal monomers to form a pathogenic--normal heterodimer, which then dissociates into two pathogenic monomers. We assume that the dissociation step is fast and  treat this conversion process as a single second-order chemical reaction, with rate $k_{12}$. Normal proteins are produced at rate $k_0$ and cleared at rate $k_1$, and pathogenic proteins are cleared at rate $\tilde{k}_1$.

The resulting reaction--diffusion system with zero-flux boundary conditions is given by:
\begin{equation}
	\label{eq:HM:strongformulation}
	\begin{cases}
		\dfrac{\partial u}{\partial t}  =\nabla \cdot (\mathbf{D} \nabla u) - k_{1} \, u - k_{12} \, u \, \tilde{u} + k_0, &  \mathrm{in} \: \Omega \times (0,T], \\[6pt]
		\dfrac{\partial \tilde{u}}{\partial t}  = \nabla \cdot (\mathbf{D} \nabla \tilde{u}) - \tilde{k}_1\, \tilde{u} + k_{12}\, \tilde{u} \, u, & \mathrm{in} \: \Omega \times (0,T], \\[6pt]
		(\mathbf{D}\nabla u )\cdot \boldsymbol{n}_\Omega = 0 \mathrm{,} \quad\;\;\, (\mathbf{D} \nabla \tilde{u} )\cdot \boldsymbol{n}_\Omega  = 0,  & \mathrm{on} \; \partial \Omega \times (0,T], \\[6pt]
		u(\boldsymbol{x},0) = u_{0}(\boldsymbol{x})\mathrm{,} \qquad \tilde{u}(\boldsymbol{x},0)  = \tilde{u}_{0}(\boldsymbol{x}), & \mathrm{in} \: \Omega,
	\end{cases}
\end{equation}
where $\mathbf{D}$ is the diffusion tensor, $\Omega$ is the spatial domain (i.e.~the brain), $\boldsymbol{n}_\Omega$ is a boundary normal, and $u_0$, $\tilde{u}_0$ are the initial protein concentrations. The model thus comprises an initial--boundary value problem for a parabolic system of two PDEs coupled by a nonlinear term.
\par
The spatial transport of A\textbeta{} proteins is a combination of diffusion in the extracellular space and axonal transport. Following \cite{weickenmeier_multiphysics_2018,weickenmeier_physics_2019}, we model this combination as anisotropic diffusion whose principal direction is aligned with the brain's axon fibers:
\begin{equation}
	\label{eq:difftensor}
	\mathbf{D}(\boldsymbol{x}) = d_\mathrm{ext}\mathbf{I} + d_\mathrm{axn}\boldsymbol{\bar{a}}(\boldsymbol{x}) \otimes \boldsymbol{\bar{a}}(\boldsymbol{x}) = d_\mathrm{ext}\mathbf{I} + d_\mathrm{axn} \mathbf{D}_\mathrm{axn}(\boldsymbol{x}).
\end{equation}
The vector field of fiber directions, $\boldsymbol{\bar{a}}=\boldsymbol{\bar{a}}(\boldsymbol{x})$, is obtained as the principal eigenvector of the diffusion tensor computed from diffusion-weighted MRI~(see \cite{corti_discontinuous_2023,mardal_mathematical_2022}). The construction of $\boldsymbol{\bar{a}}$ and choice of diffusion constants is discussed in Section~\ref{sec:num_brain}.
\par
If the kinetic rate constants are all assumed constant, then there are two spatially-homogeneous equilibria, namely:
\begin{equation}
\begin{aligned}
    &\text{healthy:} \quad &&u\equiv k_0/k_1, \quad &&\tilde{u} \equiv 0,\\
    &\text{pathogenic:} \quad &&u\equiv \tilde{k}_1/k_{12}, \quad &&\tilde{u} \equiv k_0/\tilde{k}_1 - k_1/k_{12}.
\end{aligned}
\end{equation}
Crucially, the pathogenic equilibrium is positive, and therefore physically relevant, if and only if the dimensionless \textit{basic reproduction  number} $R_0 = k_0 k_{12} / k_1 \tilde{k}_1$ is greater than unity. Indeed, as $R_0$ increases through unity, a transcritical bifurcation occurs: when $R_0 < 1$, the healthy equilibrium is stable, and no other (physical) equilibrium exists, and when $R_0 > 1$, the healthy equilibrium is unstable and the pathogenic equilibrium is stable (see, e.g., \cite{ahern_modelling_2025} for more details).
The kinetic parameters of the heterodimer formulation can be calibrated by matching model outputs to longitudinal biomarker data and established disease time scales. Such calibration can be performed either using PET measurements, as in~\cite{chaggar2025personalised,corti_uncertainty_2024} or based on employing protein concentrations derived \textit{ex vivo}, as proposed in~\cite{corti_exploring_2024}.

\subsection{Multiple-network porous media model}
\label{sec:model:MP}

For the perfusion model, our starting point is the multiple-network porous medium model proposed in \cite{corti_numerical_2023}, which is in turn inspired by \cite{jozsa_porous_2021,jozsa_sensitivity_2021,tully_cerebral_2011}.\\

We suppose there are three vascular networks, namely arterial, capillary, and venous, with corresponding pressure fields $p_\mathrm{A}(\boldsymbol{x},t)$, $p_\mathrm{C}(\boldsymbol{x},t)$, $p_\mathrm{V}(\boldsymbol{x},t)$. Because the protein concentrations $u$, $\tilde{u}$ evolve on timescales of hours to days, whereas $p_\mathrm{A}$, $p_\mathrm{C}$, $p_\mathrm{V}$ adjust to permeability changes (e.g.~due to vasoconstriction) within seconds, we assume that all pressures are quasi-static. Conservation of mass then yields the following elliptic system:
\begin{equation}
	\label{eq:MP:brainproblem}
	\begin{dcases}
		-		\nabla\cdot\left(\boldsymbol{\mathrm{K}}_\mathrm{A}\nabla p_\mathrm{A}
		\right) + \beta_\mathrm{AC}(p_\mathrm{A}-p_\mathrm{C}) = 0,       
		& \mathrm{in}\; \Omega,		
		\\[6pt]
        - \nabla\cdot\left(k_\mathrm{C} \nabla p_\mathrm{C}\right) - \beta_\mathrm{AC}(p_\mathrm{A}-p_\mathrm{C}) + \beta_\mathrm{CV}(p_\mathrm{C}-p_\mathrm{V}) = 0,
		& \mathrm{in}\; \Omega,
		\\[6pt] -\nabla\cdot\left(\mathbf{K}_\mathrm{V}\nabla p_\mathrm{V}
			\right) -
			\beta_\mathrm{CV}(p_\mathrm{C}-p_\mathrm{V})
        = 0,    
		& \mathrm{in}\; \Omega.
	\end{dcases}
\end{equation}
where $\beta_\mathrm{AC}$ and $\beta_\mathrm{CV}$ are the arterial-to-capillary and capillary-to-venous transfer coefficients. Since the brain's capillary network is isotropic, we have assumed $\boldsymbol{\mathrm{K}}_\mathrm{C} = k_\mathrm{C} \boldsymbol{\mathrm{I}}$, and we will specify the arterial and venous permeability tensors $\boldsymbol{\mathrm{K}}_\mathrm{A}$ and $\boldsymbol{\mathrm{K}}_\mathrm{V}$ on a case-by-case basis in Sections~\ref{sec:num_val} and \ref{sec:num_brain}.

Arteries spanning the pial surface of the cortex supply the brain with blood; similarly, blood is drained by the pial venous network. Accordingly, we prescribe Dirichlet boundary conditions for the arterial and venous blood pressures on the pial surface $\Gamma_\mathrm{Pial} \subset \partial \Omega$, together with a zero-flux condition for the capillary compartment:
\begin{equation}
    p_\mathrm{A} = p^\mathrm{Arteries}, \quad 
    p_\mathrm{V} = p^\mathrm{Veins}, \quad \nabla p_\mathrm{C} \cdot \boldsymbol{n}_\Omega = 0  \qquad \mathrm{on}\;\Gamma_\mathrm{Pial}.
\end{equation}
The remaining part of the boundary is the brain's ventricular surface, $\Gamma_\mathrm{Vent} = \partial \Omega \setminus \Gamma_\mathrm{Pial}$, where we prescribe zero-flux conditions:
\begin{equation}
    \nabla p_\mathrm{A} \cdot \boldsymbol{n}_\Omega = \nabla p_\mathrm{C} \cdot \boldsymbol{n}_\Omega = \nabla p_\mathrm{V} \cdot \boldsymbol{n}_\Omega = 0, \qquad \mathrm{on}\; \Gamma_\mathrm{Vent}.
\end{equation}
For the perfusion model, physiologically admissible ranges for the permeability and hemodynamic parameters can be informed by recent porous-medium perfusion studies, where these quantities are inferred directly from ASL and structural MRI data~\cite{jozsa_porous_2021,josza_mri_2023}.
\subsection{Coupled A\textbeta{} and perfusion model}
\label{sec:coupled:model}
In order to model the interaction between A\textbeta{} and the brain's microvasculature, we couple the protein and perfusion models of Sections~\ref{sec:HM:model} and \ref{sec:model:MP}, respectively. The biological mechanisms modeled here are reviewed in \cite{korte_cerebral_2020}.

First, A\textbeta{} is vasoconstrictive, i.e.~it induces blood vessels to become narrow \cite{thomas_beta_1996,Niwa2001}, thus increasing vascular resistance. In particular, Nortley et al.\ recently found that A\textbeta{} oligomers cause brain capillaries to constrict within seconds, through a mechanism involving oxidative stress and the vasoconstrictor endothelin 1 \cite{nortley_amyloid_2019}. Therefore, we assume that the permeability of the capillary bed, $k_\mathrm{C}$, is a decreasing sigmoidal function of the pathogenic protein concentration $\tilde{u}$:
\begin{equation}
    \tilde{k}_\mathrm{C}(\tilde{u}) = k_\mathrm{C}-\left(k_\mathrm{C}-k^\mathrm{A\beta}_\mathrm{C}\right)\tanh{(\alpha_{k_\mathrm{C}} \tilde{u})},
\end{equation}
where $k_\mathrm{C}$ is the base permeability (in the absence of pathogenic proteins), $k^\mathrm{A\beta}_\mathrm{C}$ is the permeability for very large $\tilde{u}$, and $\alpha_{k_\mathrm{C}}$ modulates the sensitivity to $\tilde{u}$. At the same time, we adopt a similar strategy for the transfer coefficients $\beta_\mathrm{AC}$ and $\beta_\mathrm{CV}$. In particular, for a generic couple $(i,j)$, we define:
\begin{equation}
    \tilde{\beta}_{ij}(\tilde{u}) = \beta_{ij}-\left(\beta_{ij}-\beta^\mathrm{A\beta}_{ij}\right)\tanh{(\alpha_{\beta_{ij}} \tilde{u})},
\end{equation}
where $\alpha_{\beta_{ij}}$ is a constant that modulates the impact of misfolded A\textbeta{} on the flow between compartments. Indeed, according to the multiple-network porous medium theory developed in \cite{aifantis_continuum_1979} for geophysical applications and later adapted for brain circulation in \cite{tully_cerebral_2011}, the parameter $\beta_{ij}$ models fluid exchange between networks. As discussed in \cite{peyrounette_multiscale_2018}, $\beta_{ij}$ must be proportional to the exchange area between compartments. In Alzheimer's disease, A\textbeta{} induces capillary constriction~\cite{nortley_amyloid_2019}, thus reducing capillary lumen and thus decreasing $\beta_\mathrm{AC}$ and $\beta_\mathrm{CV}$. Through these couplings, pathogenic A\textbeta{} causes local decreases in permeability (i.e.~increases in resistance), which we expect to cause hypoperfusion (reduced blood flow). Finally, the choice of  a sigmoidal function for the coupling is motivated by the experimental data  (\cite[Fig. 2D]{nortley_amyloid_2019}) that shows that soluble A\textbeta{} oligomers induce a pericyte-mediated capillary constriction with a nonlinear, saturating dependence on their concentration.

Conversely, hypoperfusion increases the production rate of A\textbeta{} \cite{shi_hypoperfusion_2000,sun_hypoxia_2006,Zhang2007}; it may also decrease its clearance rate \cite{korte_cerebral_2020}, though this is less certain. Therefore, we let the rates of A\textbeta{} production and clearance depend on the rate of blood flow, as follows. We define the ``CBF rate'' at a point $\boldsymbol{x} \in \Omega$ as the rate of flow from the arterial compartment into the capillary bed per unit mass of tissue:
\begin{equation}
    Q\big(p_\mathrm{A}(\boldsymbol{x}\big),p_\mathrm{C}(\boldsymbol{x}))=\frac{\beta_\mathrm{AC}}{\rho}\big(p_\mathrm{A}(\boldsymbol{x})-p_\mathrm{C}(\boldsymbol{x})\big),
\end{equation}
where $\rho$ is the brain tissue density estimated $1000\,\mathrm{kg\,m^{-3}}$ \cite{josza_mri_2023}.
In the framework of the multiple-network porous-medium model, this quantity coincides with the arteriole-to-capillary exchange flux per unit tissue mass, representing the effective delivery of blood to the microvascular bed at the tissue scale. {Hence}, $Q$ is a local proxy for CBF. Its spatial distribution and magnitude are determined by the pressure drop between arterial and capillary compartments and by the coupling coefficient $\beta_{\mathrm{AC}}$, which is typically chosen in literature to be consistent with perfusion values inferred from arterial spin labelling images~\cite{jozsa_porous_2021}.
We denote by $Q_\mathrm{H} = Q_\mathrm{H}(\boldsymbol{x})$ the healthy CBF rate, i.e.~corresponding to the pressure fields $p_\mathrm{A}$, $p_\mathrm{C}$ when pathogenic proteins are absent (so that there is no vasoconstriction). We now assume that the A\textbeta{} production rate increases in response to hypoperfusion, which we define as the relative decrease in $Q$ from its healthy rate $Q_\mathrm{H}$:
\begin{equation}
\label{eq:k_0^B}
    k_0^\mathrm{B}(p_\mathrm{A},p_\mathrm{C}) = k_0 + \kappa_0\left(\dfrac{Q_\mathrm{H}-Q(p_\mathrm{A},p_\mathrm{C})}{Q_\mathrm{H}}\right),
\end{equation}
where $k_0$ is the base rate, and $\kappa_0$ is the sensitivity to hypoperfusion. Similarly, we assume the clearance rates decrease in response to hypoperfusion:
\begin{equation}
    k_1^\mathrm{B}(p_\mathrm{A},p_\mathrm{C}) = k_1 - \kappa_1\left(\dfrac{Q_\mathrm{H}-Q(p_\mathrm{A},p_\mathrm{C})}{Q_\mathrm{H}}\right), \qquad \tilde{k}_1^\mathrm{B}(p_\mathrm{A},p_\mathrm{C}) = \tilde{k}_1 - \tilde{\kappa}_1\left(\dfrac{Q_\mathrm{H}-Q(p_\mathrm{A},p_\mathrm{C})}{Q_\mathrm{H}}\right).
\end{equation}
 The coupled protein--perfusion model is a parabolic--elliptic initial--boundary value problem to be solved for the concentrations $u(\boldsymbol{x},t)$, $\tilde{u}(\boldsymbol{x},t)$ and the pressures $p_\mathrm{A}(\boldsymbol{x},t)$, $p_\mathrm{C}(\boldsymbol{x},t)$, $p_\mathrm{V}(\boldsymbol{x},t)$:
\begin{equation}
	\label{eq:coupled:strongformulation}
	\begin{cases}
 	- \nabla\cdot\left(\boldsymbol{\mathrm{K}}_\mathrm{A}\nabla p_\mathrm{A}
		\right) + \tilde{\beta}_\mathrm{AC}(\tilde{u})(p_\mathrm{A}-p_\mathrm{C}) = 0       
		& \mathrm{in}\; \Omega  \times (0,T],	\\[6pt]
        - \nabla\cdot \left(\tilde{k}_\mathrm{C}(\tilde{u}) \nabla p_\mathrm{C}\right) - \tilde{\beta}_\mathrm{AC}(\tilde{u})(p_\mathrm{A}-p_\mathrm{C}) + \tilde{\beta}_\mathrm{CV}(\tilde{u})(p_\mathrm{C}-p_\mathrm{V}) = 0
		& \mathrm{in}\; \Omega \times (0,T],	\\[6pt]
        -\nabla\cdot\left(\mathbf{K}_\mathrm{V}\nabla p_\mathrm{V}\right) - \tilde{\beta}_\mathrm{CV}(\tilde{u})(p_\mathrm{C}-p_\mathrm{V})
        = 0 & \mathrm{in}\; \Omega \times (0,T], \\[6pt]
		\dfrac{\partial u}{\partial t}  =\nabla \cdot (\mathbf{D} \nabla u) - k_1^\mathrm{B}(p_\mathrm{A},p_\mathrm{C}) \, u - k_{12} \, u \, \tilde{u} + k_0^\mathrm{B}(p_\mathrm{A},p_\mathrm{C}) &  \mathrm{in} \: \Omega \times (0,T], \\[6pt]
		\dfrac{\partial \tilde{u}}{\partial t}  = \nabla \cdot (\mathbf{D} \nabla \tilde{u}) - \tilde{k}_1^\mathrm{B}(p_\mathrm{A},p_\mathrm{C})\, \tilde{u} + k_{12}\, \tilde{u} \, u & \mathrm{in} \: \Omega \times (0,T].\\[6pt]
    \end{cases}
\end{equation}
The system is complemented with the following boundary conditions:
\begin{equation}
\begin{cases}
        p_\mathrm{A} = p^\mathrm{Arteries}, \qquad\;\;\; 
        p_\mathrm{V} = p^\mathrm{Veins}, & \mathrm{on}\;\Gamma_\mathrm{Pial}\times (0,T],\\[6pt]
        \nabla p_\mathrm{A} \cdot \boldsymbol{n}_\Omega = \nabla p_\mathrm{V} \cdot \boldsymbol{n}_\Omega = 0 & \mathrm{on}\; \Gamma_\mathrm{Vent}\times (0,T],\\[6pt]
		\nabla p_\mathrm{C} \cdot \boldsymbol{n}_\Omega =0, \qquad (\mathbf{D}\nabla u )\cdot \boldsymbol{n}_\Omega = 0 \mathrm{,} \qquad\; (\mathbf{D} \nabla \tilde{u} )\cdot \boldsymbol{n}_\Omega  = 0  & \mathrm{on} \; \partial \Omega \times (0,T],
    \end{cases}
\end{equation}
and with the initial conditions
\begin{equation}
\label{eq:coupled:init_cond}
    u(0,\boldsymbol{x}) = u_{0}(\boldsymbol{x})\mathrm{,} \qquad \tilde{u}(0,\boldsymbol{x})  = \tilde{u}_{0}(\boldsymbol{x}) \quad \mathrm{in} \: \Omega.
\end{equation}

\subsubsection*{Nondimensionalisation}
For the heterodimer model, we rescale the variables as follows, where tildes indicate the new dimensionless variables \cite{ahern_modelling_2025}:
\begin{equation*}
    u = \frac{k_0}{k_1} \hat{u}, \qquad \tilde{u} = \frac{k_0}{\tilde{k}_1} \hat{\tilde{u}}, \qquad t = \frac{1}{\tilde{k}_1} \hat{t}, \qquad x = \sqrt{\frac{d_\mathrm{ext}}{\tilde{k}_1}} \hat{x}.
\end{equation*}
Moreover, we scale the pressure fields so as to map the interval $[p^\mathrm{Veins},p^\mathrm{Arteries}]$  uniformly onto $[0,1]$:
\begin{equation*}
    p_\mathrm{A} = (p^\mathrm{Arteries}-p^\mathrm{Veins}) \hat{p}_\mathrm{A} + p^\mathrm{Veins}, \quad p_\mathrm{C} = (p^\mathrm{Arteries}-p^\mathrm{Veins}) \hat{p}_\mathrm{C} + p^\mathrm{Veins}, \quad p_\mathrm{V} = (p^\mathrm{Arteries}-p^\mathrm{Veins}) \hat{p}_\mathrm{V} + p^\mathrm{Veins}.
\end{equation*}
Rescaling the equations of the system \eqref{eq:coupled:strongformulation}, we obtain the dimensionless system:
\begin{equation}
	\label{eq:coupled:dimlessstrongformulation}
	\begin{cases}
 	- \nabla\cdot\left(\boldsymbol{\sigma}_\mathrm{A}\nabla p_\mathrm{A}
		\right) + \tilde{\gamma}_\mathrm{AC}(\tilde{u})(p_\mathrm{A}-p_\mathrm{C}) = 0       
		& \mathrm{in}\; \Omega  \times (0,T],	\\[6pt]
        - \nabla\cdot \left(\tilde{\sigma}_\mathrm{C}(\tilde{u}) \nabla p_\mathrm{C}\right) - \tilde{\gamma}_\mathrm{AC}(\tilde{u})(p_\mathrm{A}-p_\mathrm{C}) + B\tilde{\gamma}_\mathrm{CV}(\tilde{u})(p_\mathrm{C}-p_\mathrm{V}) = 0
		& \mathrm{in}\; \Omega \times (0,T],	\\[6pt]
        -\nabla\cdot\left(\boldsymbol{\sigma}_\mathrm{V}\nabla p_\mathrm{V}\right) - B\tilde{\gamma}_\mathrm{CV}(\tilde{u})(p_\mathrm{C}-p_\mathrm{V})
        = 0 & \mathrm{in}\; \Omega \times (0,T], \\[6pt]
		\epsilon\dfrac{\partial u}{\partial t}  =\nabla \cdot (\epsilon\boldsymbol{\delta} \nabla u) -\lambda_1^\mathrm{B}(p_\mathrm{A},p_\mathrm{C}) u - R \, u \, \tilde{u} + \mu_0^\mathrm{B}(p_\mathrm{A},p_\mathrm{C}) &  \mathrm{in} \: \Omega \times (0,T], \\[6pt]
		\;\,\dfrac{\partial \tilde{u}}{\partial t}  = \nabla \cdot (\boldsymbol{\delta} \nabla \tilde{u}) - \tilde{\lambda}_1^\mathrm{B}(p_\mathrm{A},p_\mathrm{C}) \tilde{u} + R\, \tilde{u} \, u & \mathrm{in} \: \Omega \times (0,T],
	\end{cases}
\end{equation}
where
\begin{equation*}
\begin{aligned}
    \boldsymbol{\sigma}_\mathrm{A} &= \dfrac{\tilde{k}_1\boldsymbol{\mathrm{K}}_\mathrm{A}}{d_\mathrm{ext} \beta_\mathrm{AC}}, 
    & B &= \dfrac{\beta_\mathrm{CV}}{\beta_\mathrm{AC}}, 
    & \epsilon &= \dfrac{\tilde{k}_1}{k_1}, 
    & \lambda_1^\mathrm{B}(p_\mathrm{A},p_\mathrm{C}) &= \dfrac{1}{k_1} k_1^\mathrm{B}(p_\mathrm{A},p_\mathrm{C}), 
    \\
    \tilde{\sigma}_\mathrm{C}(\tilde{u}) &= \dfrac{\tilde{k}_1\,\tilde{k}_\mathrm{C}(\tilde{u})}{d_\mathrm{ext} \beta_\mathrm{AC}}, 
    & \tilde{\gamma}_\mathrm{AC}(\tilde{u}) &= \dfrac{1}{\beta_{\mathrm{AC}}}\tilde{\beta}_\mathrm{AC}(\tilde{u}), 
    & R_0 &= \dfrac{k_0 k_{12}}{k_1 \tilde{k}_1}, 
    & \tilde{\lambda}_1^\mathrm{B}(p_\mathrm{A},p_\mathrm{C}) &= \dfrac{1}{\tilde{k}_1} \tilde{k}_1^\mathrm{B}(p_\mathrm{A},p_\mathrm{C}), 
    \\
    \boldsymbol{\sigma}_\mathrm{V} &= \dfrac{\tilde{k}_1\boldsymbol{\mathrm{K}}_\mathrm{V}}{d_\mathrm{ext} \beta_\mathrm{AC}}, 
    & \tilde{\gamma}_\mathrm{CV}(\tilde{u}) &= \dfrac{1}{\beta_{\mathrm{CV}}}\tilde{\beta}_\mathrm{CV}(\tilde{u}), 
    & \boldsymbol{\delta} &= \dfrac{\mathbf{D}}{d_\mathrm{ext}}, 
    & \mu_0^\mathrm{B}(p_\mathrm{A},p_\mathrm{C}) &= \dfrac{1}{k_0} k_0^\mathrm{B}(p_\mathrm{A},p_\mathrm{C}).
\end{aligned}
\end{equation*}
The system is complemented with the following boundary conditions:
\begin{equation}
	\label{eq:coupled:dimlessbcs}
	\begin{cases}     
        p_\mathrm{A} = 1, \qquad 
        p_\mathrm{V} = 0, & \mathrm{on}\;\Gamma_\mathrm{Pial}\times (0,T],\\[6pt]
        \nabla p_\mathrm{A} \cdot \boldsymbol{n} = \nabla p_\mathrm{V} \cdot \boldsymbol{n} = 0, & \mathrm{on}\; \Gamma_\mathrm{Vent}\times (0,T],\\[6pt]
        \nabla p_\mathrm{C} \cdot \boldsymbol{n} = 0, \qquad
		(\epsilon\boldsymbol{\delta}\nabla u )\cdot \boldsymbol{n} = 0 \mathrm{,} \qquad (\boldsymbol{\delta} \nabla \tilde{u} )\cdot \boldsymbol{n}  = 0,  & \mathrm{on} \; \partial \Omega \times (0,T],
	\end{cases}
\end{equation}
and the initial data $u_0$ and $\tilde{u}_{0}$ in Equation \eqref{eq:coupled:init_cond} have been suitably rescaled and nondimensionalized.
\begin{remark}
In  \eqref{eq:coupled:dimlessstrongformulation}, the parameter $\epsilon$ represents the ratio between the characteristic clearance time of the pathogenic species and that of the healthy one. In line with the analysis of the heterodimer--vascular model in \cite{ahern_modelling_2025}, we typically have $\epsilon = O(1)$ and $\epsilon$ does not act as a small parameter. This choice reflects biologically informed estimates of the underlying kinetic rates (see also the sensitivity analysis in \cite{corti_exploring_2024}). Consequently, we do not expect a strong separation of time scales between healthy and pathogenic protein dynamics (which would lead to singular perturbation limit), and the {stability analysis} reported in this work arises from the nonlinear coupling mechanisms.
\end{remark}
\section{Discrete formulation}
\label{sec:discretization}
In this section, we introduce the discretization scheme for the resolution of the system in equation~\eqref{eq:coupled:strongformulation}. In particular, in Section~\ref{sec:space_discretization}, we discuss the discretization in space, by means of a dG (discontinuous Galerkin) method, and in Section~\ref{sec:time_discretization}, we discretize in time by means of an implicit Euler time stepping algorithm.

\subsection{Space discretization: discontinuous Galerkin method}
\label{sec:space_discretization}
We first introduce a mesh partition $\partition$ of the domain $\Omega$ made of shape-regular simplicial elements $K$, with corresponding diameter $h_K$ and define $ h=\max_{K\in\partition}\{h_K\}$. If two elements share a common face $F$ (e.g.,~the triangular interface of two tetrahedra), we call it an interior face, $F\in\facesinternal$. Otherwise, we call $F$ a boundary face, $F\in\facesboundary$. The set is partitioned on the pial and ventricular regions of $\partial \Omega$, i.e.~$\facesboundary = \facesboundaryvent \cup \facesboundarypial$. We now introduce the so-called trace operators. Let $F\in\facesinternal$ be a face shared by two elements $K^\pm$. We denote by $\boldsymbol{n}^\pm$ the unit normal vector to $F$ pointing outward to $K^\pm$, respectively. Then, for sufficiently regular scalar-valued $q$ and vector-valued functions $\boldsymbol{v}$, respectively, we define the average operator $\averagel{\cdot}\averager$ as $\averagel{q}\averager = (q^+ + q^-)/2$, and $\averagel{\boldsymbol{v}}\averager = (\boldsymbol{v}^+ + \boldsymbol{v}^-)/2$, and the jump operator $\jumpl{\cdot}\jumpr$ as  $\jumpl{q}\jumpr = q^+\boldsymbol{n}^+ + q^-\boldsymbol{n}^-$, and $\jumpl{\boldsymbol{v}}\jumpr = \boldsymbol{v}^+\cdot\boldsymbol{n}^+ + \boldsymbol{v}^-\cdot\boldsymbol{n}^-$. The $\pm$ superscripts denote the traces on the face $F$ of the functions defined on $K^\pm$. 
Analogously, on the face $F\in\facesboundary$ of a cell $K\in\partition$, we define the average operator $\averagel{\cdot}\averager$ as $\averagel{q}\averager = q$ and $\averagel{\boldsymbol{v}}\averager = \boldsymbol{v}$, and the jump operator $\jumpl{\cdot}\jumpr$ as $\jumpl{q}\jumpr = (q-g)\boldsymbol{n}$ and $\jumpl{\boldsymbol{v}}\jumpr = (\boldsymbol{v}-\boldsymbol{g})\cdot\boldsymbol{n}$, where 
$g$ and $\boldsymbol{g}$ are regular enough Dirichlet boundary data and  $\boldsymbol{n}$
$\boldsymbol{n}$ is the outward unit normal vector to $\partial\Omega$.
\par
Let us define $\mathbb{P}_{\ell}(K)$ as the space of polynomials of total degree $\ell\geq 1$ over a mesh element $K$. Then we can introduce the following definitions of dG finite element spaces $\Vh = \{v\in L^2(\Omega):\; v|_K\in\mathbb{P}_{\ell}(K)\;\forall K\in\partition\}$. From now on, we introduce also a shorthand notation for the integrals over the faces $\int_\mathcal{F} = \sum_{F\in\mathcal{F}}\int_F$, and the notation $\nabla_h$ to denote the broken gradient operator. Next, we introduce the following forms for all $v,w, \tilde{u}\in\Vh$:
\begin{align*}
    \mathcal{A}_j(v,w) = & \int_{\Omega}\mathbf{K}_j\nabla_h v\cdot\nabla_h w +\int_{\facesinternal\cup\facesboundarypial}\left(\eta_j \jumpl v\jumpr \cdot \jumpl w\jumpr -\averagel\mathbf{K}_j \nabla_{h} v\averager \cdot \jumpl w \jumpr -  \jumpl v\jumpr \cdot \averagel\mathbf{K}_j \nabla_{h} w\averager\right)\mathrm{d}\sigma, & j=\mathrm{A,V} \\
    \mathcal{A}_\mathrm{C}(v,w;\tilde{u}) = & \int_{\Omega}k_\mathrm{c}(\tilde{u})\nabla_h v\cdot\nabla_h w +\int_{\facesinternal}\left(\eta_\mathrm{C}(\tilde{u}) \jumpl v\jumpr \cdot \jumpl w\jumpr -\averagel k_\mathrm{c}(\tilde{u}) \nabla_{h} v\averager \cdot \jumpl w \jumpr -  \jumpl v\jumpr \cdot \averagel k_\mathrm{c}(\tilde{u}) \nabla_{h} w\averager\right)\mathrm{d}\sigma, \\
    \mathcal{A}_\mathrm{H}(v,w) = &\int_{\Omega}\mathbf{D}\nabla_h v\cdot\nabla_h w +\int_{\facesinternal}\left(\eta_\mathrm{H} \jumpl v\jumpr \cdot \jumpl w\jumpr -\averagel\mathbf{D} \nabla_{h} v\averager \cdot \jumpl w \jumpr -  \jumpl v\jumpr \cdot \averagel\mathbf{D} \nabla_{h} w\averager\right)\mathrm{d}\sigma,
\end{align*}
where $\eta_j$ with $j=\mathrm{A,C,V,H}$ are the discontinuity penalization functions defined as follows:
\begin{align}
\label{penltycoeff}
  \eta_\mathrm{A} = & \eta_{0\mathrm{A}} 
  \begin{dcases}
    \{\|\mathbf{K}_\mathrm{A}\|\}_\mathrm{harm}\dfrac{\ell^2}{\{h_K\}_\mathrm{harm}} & \mathrm{on} \; F \in \facesinternal, \\
    \|\mathbf{K}_\mathrm{A}\|\dfrac{\ell^2}{h_K} & \mathrm{on} \; F \in \facesboundarypial,
  \end{dcases} \quad && \eta_\mathrm{V} = \eta_{0\mathrm{V}} 
  \begin{dcases}
    \{\|\mathbf{K}_\mathrm{V}\|\}_\mathrm{harm}\dfrac{\ell^2}{\{h_K\}_\mathrm{harm}} & \mathrm{on} \; F \in \facesinternal, \\
    \|\mathbf{K}_\mathrm{V}\|\dfrac{\ell^2}{h_K} & \mathrm{on} \; F \in \facesboundarypial,
  \end{dcases} \\
    \eta_\mathrm{C}(\tilde{u}) = & \eta_{0\mathrm{C}} 
    \{k_\mathrm{c}(\tilde{u})\}_\mathrm{harm}\dfrac{\ell^2}{\{h_K\}_\mathrm{harm}} \qquad \mathrm{on} \; F \in \facesinternal &&
    \eta_\mathrm{H} = \eta_{0\mathrm{H}} 
    \{\|\mathbf{D}\|\}_\mathrm{harm}\dfrac{\ell^2}{\{h_K\}_\mathrm{harm}} \;\qquad \mathrm{on} \; F \in \facesinternal.
\end{align}
where $\eta_{0j}$  with $j=\mathrm{A,C,V,H}$ are constant parameters that should be chosen sufficiently large to ensure the stability of the discrete formulation (see \cite{corti_numerical_2023,antonietti_discontinuous_2024} for details), $\|\cdot\|$ is the euclidean matrix norm, and $\{\cdot\}_\mathrm{harm}$ is the harmonic mean.
The semi-discrete formulation in space reads as follows. Given the initial conditions $u_h(0) = u_{0h}\in\Vh$ and $\tilde{u}_h(0)  = \tilde{u}_{0h}\in\Vh$, for each $t \in (0,T)$, find $(p_{\mathrm{A}h}(t),p_{\mathrm{C}h}(t),p_{\mathrm{V}h}(t), u_h(t), \tilde{u}_h(t))\in \Vh\times\Vh\times\Vh\times\Vh\times\Vh$ such that:
\begin{subequations}
\label{eq:coupled:semidicreteformulation}
\begin{alignat}{3}
    \mathcal{A}_\mathrm{A}(p_{\mathrm{A}h},q_{\mathrm{A}h}) + (\tilde{\beta}_\mathrm{AC}(\tilde{u}_h)(p_{\mathrm{A}h}-p_{\mathrm{C}h}),q_{\mathrm{A}h})_\Omega = & \; 0,    
    & & \quad \forall q_\mathrm{A}\in \Vh,	\\[6pt]
    \mathcal{A}_\mathrm{C}(p_{\mathrm{C}h},q_{\mathrm{C}h};\tilde{u}_h)  - (\tilde{\beta}_\mathrm{AC}(\tilde{u}_h)(p_{\mathrm{A}h}-p_{\mathrm{C}h}),q_{\mathrm{C}h})_\Omega + (\tilde{\beta}_\mathrm{CV}(\tilde{u}_h)(p_{\mathrm{C}h}-p_{\mathrm{V}h}),q_{\mathrm{C}h})_\Omega = & \, 0,
	& & \quad \forall q_\mathrm{C}\in \Vh,	\\[6pt]
    \mathcal{A}_\mathrm{V}(p_{\mathrm{V}h},q_{\mathrm{V}h}) - (\tilde{\beta}_\mathrm{CV}(\tilde{u}_h)(p_{\mathrm{C}h}-p_{\mathrm{V}h}),q_{\mathrm{V}h})_\Omega = & \, 0, & & \quad  \forall q_\mathrm{V}\in \Vh, \\[6pt]
	\left(\dot{u}_h,v_h\right)_\Omega +\mathcal{A}_\mathrm{H}(u_h, v_h) + (k_1^\mathrm{B}(p_{\mathrm{A}h},p_{\mathrm{C}h})u_h-k_0^\mathrm{B}(p_{\mathrm{A}h},p_{\mathrm{C}h}),v_h)_\Omega + \left(k_{12} \, \tilde{u}_h u_h, v_h\right)_\Omega = &  \, 0, & & \quad \forall v\in \Vh, \\[6pt]
	   \left(\dot{\tilde{u}}_h,\tilde{v}_h\right)_\Omega+ \mathcal{A}_\mathrm{H}(\tilde{u}_h, \tilde{v}_h) + (\tilde{k}_1^\mathrm{B}(p_{\mathrm{A}h},p_{\mathrm{C}h})\, \tilde{u}_h,\tilde{v}_h)_\Omega - \left(k_{12} u_h\, \tilde{u}_h,\tilde{v}_h\right)_\Omega = & \, 0, & & \quad \forall \tilde{v}\in \Vh.
\end{alignat}
\end{subequations}

\subsection{Time discretization: implicit Euler finite difference scheme}
\label{sec:time_discretization}
Before discretizing in time the problem in Equation \eqref{eq:coupled:semidicreteformulation}, we construct the matrices associated with the problem.
Let $\{\varphi_{j}\}^{N_h}_{j=1}$ be suitable basis functions for the discrete spaces $\Vh$. Then we can write the unknowns as an expansion in the polynomial basis:
\begin{equation*}
   u_h(\mathbf{x},t) = \sum_{j=1}^{N_h} U_j(t)\varphi_j(\mathbf{x}), \qquad 
   \tilde{u}_h(\mathbf{x},t) = \sum_{j=1}^{N_h} \tilde{U}_j(t)\varphi_j(\mathbf{x}), \qquad
   p_{kh}(\mathbf{x},t) = \sum_{j=1}^{N_h} P_{kj}(t)\varphi_j(\mathbf{x}), \qquad k=\{\mathrm{A,C,V}\}.
\end{equation*}
We denote by $\mathbf{U} = [U_j]_{j=1}^{N_h} \in \mathbb{R}^{N_h}$, $\widetilde{\mathbf{U}} = [\tilde{U}_j]_{j=1}^{N_h} \in \mathbb{R}^{N_h}$, $\mathbf{P}_k = [P_{kj}]_{j=1}^{N_h} \in \mathbb{R}^{N_h}$ the vectors of the expansion coefficients, and define the following matrices for $i,j = 1,...,N_h$:
\begin{subequations}
\begin{alignat}{2}
    &[\mathbf{A}_\mathrm{A}]_{ij} =  \mathscr{A}_\mathrm{A}(\varphi_j, \varphi_i) & (\text{Arterial stiffness matrix}); \\
    &[\mathbf{A}_\mathrm{C}(\widetilde{\mathbf{U}}_{h})]_{ij} =  \mathscr{A}_\mathrm{C}(\varphi_j, \varphi_i; \tilde{u}_h) & (\text{Capillary stiffness matrix}); \\
    &[\mathbf{A}_\mathrm{V}]_{ij} =  \mathscr{A}_\mathrm{V}(\varphi_j, \varphi_i) & (\text{Venous stiffness matrix}); \\
    &[\mathbf{M}_{\beta_\mathrm{AC}}(\widetilde{\mathbf{U}}_{h})]_{ij} =  (\tilde{\beta}_\mathrm{AC}(\tilde{u}_h)\varphi_j, \varphi_i)_\Omega & (\text{Arterial--capillary\;coupling\;matrix}); \\
    &[\mathbf{M}_{\beta_\mathrm{CV}}(\widetilde{\mathbf{U}}_{h})]_{ij} =  (\tilde{\beta}_\mathrm{CV}(\tilde{u}_h)\varphi_j, \varphi_i)_\Omega & (\text{Capillary--venous\;coupling\;matrix}); \\
    &[\mathbf{M}]_{ij} =  (\varphi_j, \varphi_i)_\Omega & (\mathrm{Mass\;matrix}); \\
    &[\mathbf{A}_\mathrm{H}]_{ij} =  \mathscr{A}_\mathrm{H}(\varphi_j, \varphi_i) & (\text{Heterodimer stiffness matrix}); \\
    &[\mathbf{M}_{k_1}(\mathbf{P}_{\mathrm{A}h},\mathbf{P}_{\mathrm{C}h})]_{ij} =  ({k}_1^\mathrm{B}(p_{\mathrm{A}h},p_{\mathrm{C}h})\, \varphi_j, \varphi_i)_\Omega & (\text{Healthy protein clearance matrix}); \\
    &[\mathbf{M}_{\tilde{k}_1}(\mathbf{P}_{\mathrm{A}h},\mathbf{P}_{\mathrm{C}h})]_{ij} =  (\tilde{k}_1^\mathrm{B}(p_{\mathrm{A}h},p_{\mathrm{C}h})\, \varphi_j, \varphi_i)_\Omega\quad & (\text{Misfolded protein clearance matrix}); \\
    &[\mathbf{M}_{k_{12}}(\mathbf{V}_{h})]_{ij} =  (k_{12}v_h\, \varphi_j, \varphi_i)_\Omega & (\text{Protein misfolding matrix}).
\end{alignat}
\end{subequations}
Moreover we define the forcing term $\mathbf{F}_{k_0}(\mathbf{P}_{\mathrm{A}h},\mathbf{P}_{\mathrm{C}h}) = [(k_0^\mathrm{B}(p_{\mathrm{A}h},p_{\mathrm{C}h}),\varphi_j)_\Omega]_{j=1}^{N_h}$. 
Finally, we adopt an implicit Euler scheme to discretize the problem in time. We construct a partition of the interval $[0, T]$ by defining $N_T$ intervals $0=t_0<t_1<...<t_{N_T}=T$. We assume a constant timestep $\Delta t = t_{n+1}-t_n$, $n=0,...,N_T-1$. Moreover, we split the resolution of the nonlinear problem into two steps, using a loosely-coupled scheme for the resolution of the porous media problem and the heterodimer one. Given $\mathbf{U}_{h}^n$ and $\widetilde{\mathbf{U}}_h^n$, solve
\begin{subequations} 
\label{eq:coupled:fullydiscreteformulation}
	\begin{alignat}{2}
    \mathbf{A}_\mathrm{A}\mathbf{P}_{\mathrm{A}h}^{n+1} 
    + \mathbf{M}_{\beta_\mathrm{AC}}(\widetilde{\mathbf{U}}_{h}^*)\left(\mathbf{P}_{\mathrm{A}h}^{n+1} -\mathbf{P}_{\mathrm{C}h}^{n+1}\right) = & \, \mathbf{0},      
    \\[6pt]
   \mathbf{A}_\mathrm{C}(\widetilde{\mathbf{U}}_{h}^*)\mathbf{P}_{\mathrm{C}h}^{n+1} 
   - \mathbf{M}_{\beta_\mathrm{AC}}(\widetilde{\mathbf{U}}_{h}^*)\left(\mathbf{P}_{\mathrm{A}h}^{n+1}  -\mathbf{P}_{\mathrm{C}h}^{n+1} \right)
   + \mathbf{M}_{\beta_\mathrm{CV}}(\widetilde{\mathbf{U}}_{h}^*)\left(\mathbf{P}_{\mathrm{C}h}^{n+1}  -\mathbf{P}_{\mathrm{V}h}^{n+1} \right) = & \, \mathbf{0},       
    \\[6pt]
   \mathbf{A}_\mathrm{V}\mathbf{P}_{\mathrm{V}h}^{n+1}  
   - \mathbf{M}_{\beta_\mathrm{CV}}(\widetilde{\mathbf{U}}_{h}^*)\left(\mathbf{P}_{\mathrm{C}h}^{n+1}  -\mathbf{P}_{\mathrm{V}h}^{n+1} \right) = & \, \mathbf{0},
   \\[6pt]
	  \left(\mathbf{M}  
       + \Delta t \left(\mathbf{A}_\mathrm{H} 
       + \mathbf{M}_{k_1}(\mathbf{P}_{\mathrm{A}h}^{n+1},\mathbf{P}_{\mathrm{C}h}^{n+1}) 
       + \mathbf{M}_{k_{12}}(\widetilde{\mathbf{U}}_{h}^*)\right)\right)\mathbf{U}_h^{n+1}
       = & \,\mathbf{M}{\mathbf{U}}_h^{n} + \Delta t \mathbf{F}_{k_0}(\mathbf{P}_{\mathrm{A}h}^{n+1},\mathbf{P}_{\mathrm{C}h}^{n+1}),
    \\[6pt]
	  \left(\mathbf{M}  
       + \Delta t \left(\mathbf{A}_\mathrm{H} 
       + \mathbf{M}_{\tilde{k}_1}(\mathbf{P}_{\mathrm{A}h}^{n+1},\mathbf{P}_{\mathrm{C}h}^{n+1}) 
       - \mathbf{M}_{k_{12}}(\widetilde{\mathbf{U}}_{h}^*)\right)\right)\widetilde{\mathbf{U}}_h^{n+1}
       = & \, \mathbf{M}{\widetilde{\mathbf{U}}}_h^{n}.
	\end{alignat}
\end{subequations}
In practice, we adopt a \emph{loosely coupled} splitting strategy by choosing $\widetilde{\mathbf{U}}_{h}^* = \widetilde{\mathbf{U}}_{h}^n$. This is equivalent to solving the porous-medium problem at time $t^{n+1}$ by using the pathogenic concentration $\widetilde{\mathbf{U}}_{h}^n$ at time $t^n$:
\begin{equation*} 
\begin{bmatrix}
    \mathbf{A}_\mathrm{A} + \mathbf{M}_{\beta_\mathrm{AC}}(\widetilde{\mathbf{U}}_{h}^n)
    & - \mathbf{M}_{\beta_\mathrm{AC}}(\widetilde{\mathbf{U}}_{h}^n)
    & \mathbf{0}      
    \\[6pt]
    - \mathbf{M}_{\beta_\mathrm{AC}}(\widetilde{\mathbf{U}}_{h}^n)
    & \mathbf{A}_\mathrm{C}(\widetilde{\mathbf{U}}_{h}^n)+\mathbf{M}_{\beta_\mathrm{AC}}(\widetilde{\mathbf{U}}_{h}^n) + \mathbf{M}_{\beta_\mathrm{CV}}(\widetilde{\mathbf{U}}_{h}^n)
    & -\mathbf{M}_{\beta_\mathrm{CV}}(\widetilde{\mathbf{U}}_{h}^n)
    \\[6pt]
    \mathbf{0}
    & - \mathbf{M}_{\beta_\mathrm{CV}}(\widetilde{\mathbf{U}}_{h}^n)
    & \mathbf{A}_\mathrm{V} + \mathbf{M}_{\beta_\mathrm{CV}}(\widetilde{\mathbf{U}}_{h}^n)
\end{bmatrix}
\begin{bmatrix}
    \mathbf{P}_{\mathrm{A}h}^{n+1} \\[6pt]
    \mathbf{P}_{\mathrm{C}h}^{n+1} \\[6pt]
    \mathbf{P}_{\mathrm{V}h}^{n+1}
\end{bmatrix} = 
\begin{bmatrix}
    \mathbf{0} \\[6pt]
    \mathbf{0} \\[6pt]
    \mathbf{0}
\end{bmatrix}.
\end{equation*}
Then, we solve the heterodimer model at the same time level using the updated pressures $p_\mathrm{A}^{n+1}$ and $p_\mathrm{C}^{n+1}$, with a decoupled strategy for the two variables (see \cite{antonietti_discontinuous_2024}):
\begin{equation*} 
\begin{aligned}
    \left(\mathbf{M}  
       + \Delta t \left(\mathbf{A}_\mathrm{H} 
       + \mathbf{M}_{k_1}(\mathbf{P}_{\mathrm{A}h}^{n+1},\mathbf{P}_{\mathrm{C}h}^{n+1}) 
       + \mathbf{M}_{k_{12}}(\widetilde{\mathbf{U}}_{h}^n)\right)\right)\mathbf{U}_h^{n+1} 
    & = \mathbf{M}\mathbf{U}_h^{n} + \Delta t \,\mathbf{F}_{k_0}(\mathbf{P}_{\mathrm{A}h}^{n+1},\mathbf{P}_{\mathrm{C}h}^{n+1}),
    \\[6pt]
    \left(\mathbf{M}  
       + \Delta t \left(\mathbf{A}_\mathrm{H} 
       + \mathbf{M}_{\tilde{k}_1}(\mathbf{P}_{\mathrm{A}h}^{n+1},\mathbf{P}_{\mathrm{C}h}^{n+1}) 
       - \mathbf{M}_{k_{12}}(\mathbf{U}_{h}^n)\right)\right)\widetilde{\mathbf{U}}_h^{n+1} 
    & = \mathbf{M}\widetilde{\mathbf{U}}_h^{n}.
\end{aligned}       
\end{equation*}
Another possible choice, although less convenient from a computational point of view, would be to take $\widetilde{\mathbf{U}}_{h}^* = \widetilde{\mathbf{U}}_{h}^{n+1}$ and linearize the resulting problem within an iterative nonlinear solver in a fully monolithic setting.
\section{Numerical validation in simple geometries}
\label{sec:num_val}
We now present some simple test cases to highlight some properties of the model introduced in Equation~\eqref{eq:coupled:strongformulation}. We consider an idealized rectangular domain $\Omega=(0\,\mathrm{m},0.1\,\mathrm{m})\times(0\,\mathrm{m},0.4\,\mathrm{m})$ and a final time $T=100\,\mathrm{years}$. We impose Dirichlet boundary conditions on $\partial\Omega$ for the pressures, representing constant arterial and venous blood pressures at the brain--body vascular interface. We use a structured triangular mesh with $15\,000$ elements, a polynomial degree $\ell=1$ for all five discretized fields, and a time step of $\Delta t=0.05$ years. Unless otherwise indicated, all spatial coordinates in this section represent values in metres ($\mathrm{m}$). 
\subsection{Test Case 1: Multistability dependent on the initial condition}
In this test case, we consider a rectangular domain and prescribe initial data with a localized region of high pathogenic protein concentration, called the ``seeding region''. We are interested in whether all such ``pathogenic seeds'' trigger disease spread, as is the case for most reaction--diffusion models of A\textbeta{} in the literature \cite{weickenmeier_physics_2019,fornari_prion-like_2019,antonietti_discontinuous_2024,corti_exploring_2024}. The parameter values used in this test case are given in Appendix A (Table \ref{tab:phys}).

\begin{figure}[t!]
	\centering
	{\includegraphics[width=\textwidth]{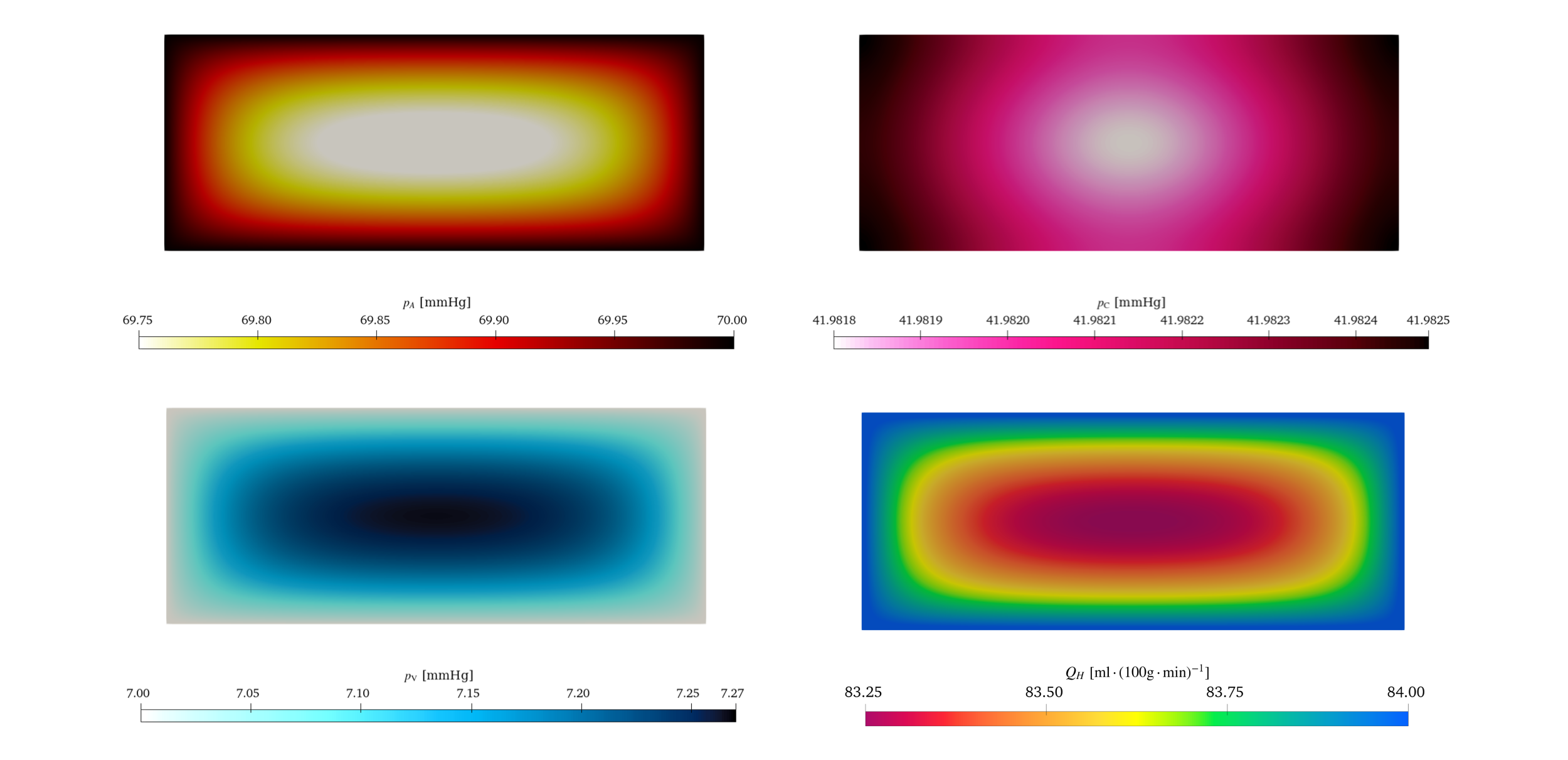}}
	\caption{Test Case 1: Pressures computed in healthy conditions in the domain $\Omega$. We report the arterial $p_\mathrm{A}$ (upper-left), capillary $p_\mathrm{C}$ (upper-right), and venous $p_\mathrm{V}$ (lower-left) pressures and the healthy CBF rate $Q_\mathrm{H}=\frac{\beta_\mathrm{AC}}{\rho}(p_\mathrm{A}-p_\mathrm{C})$ (lower right).}
	\label{fig:multiple_healthypressures}
\end{figure}

We must first compute the healthy perfusion rate $Q_\mathrm{H}$ from Equation~\eqref{eq:MP:brainproblem}, which describes the heathy CBF rate in the absence of A\textbeta{} or other vascular pathologies. Figure~\ref{fig:multiple_healthypressures} shows the numerical solution. The gradients in arterial and venous pressure near the boundary determine the flux into and out of the domain. The capillary pressure is spatially almost constant because the the value used for the permeability $k_\text{C}$ (reported in Table~\ref{tab:phys} in Appendix~\ref{app:param_values}) is relatively large.
\begin{figure}[t!]
	\centering
	{\includegraphics[width=\textwidth]{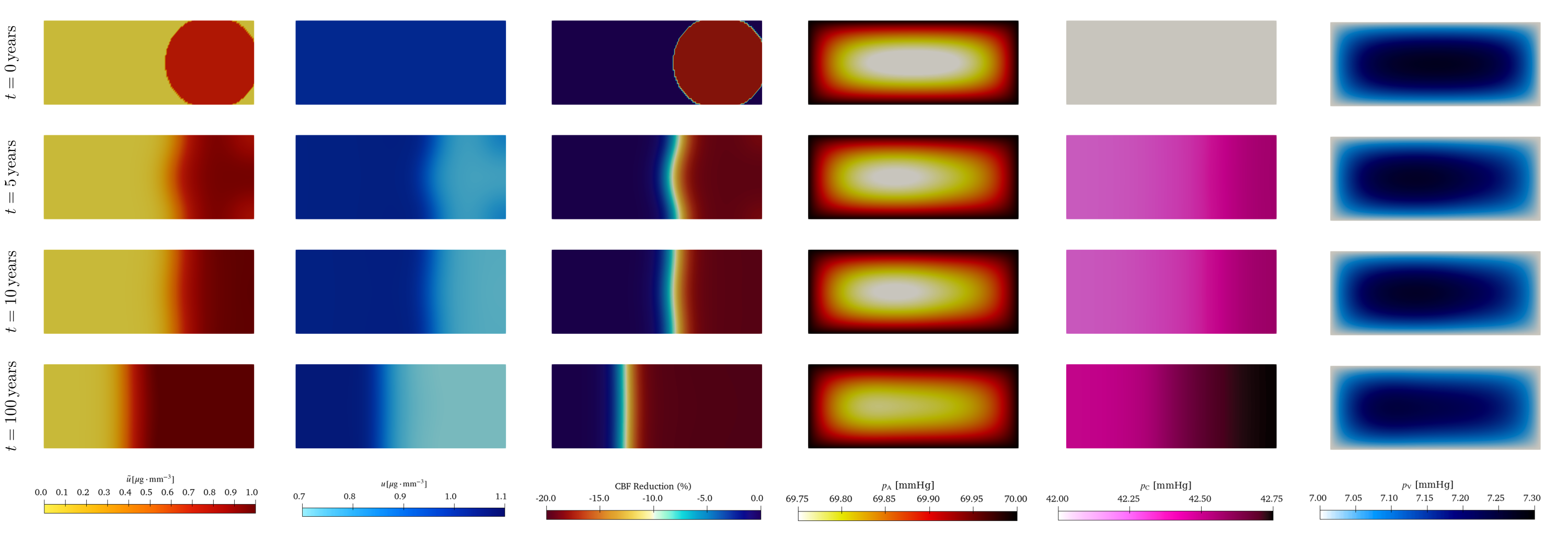}}
	\caption{Numerical solution for Test Case 1 with a large seeding region. From left to right we have misfolded proteins $\tilde{u}$, healthy proteins $u$, reduction of CBF, arterial pressure $p_\mathrm{A}$, capillary pressure $p_\mathrm{C}$, and venous pressure $p_\mathrm{V}$. Disease propagation \textit{succeeds}.  Values mapped to white in the color scale are shown in gray to improve background contrast.}
	\label{fig:multiple_wavepropagation}
\end{figure}
\begin{figure}[t!]
	\centering
	{\includegraphics[width=\textwidth]{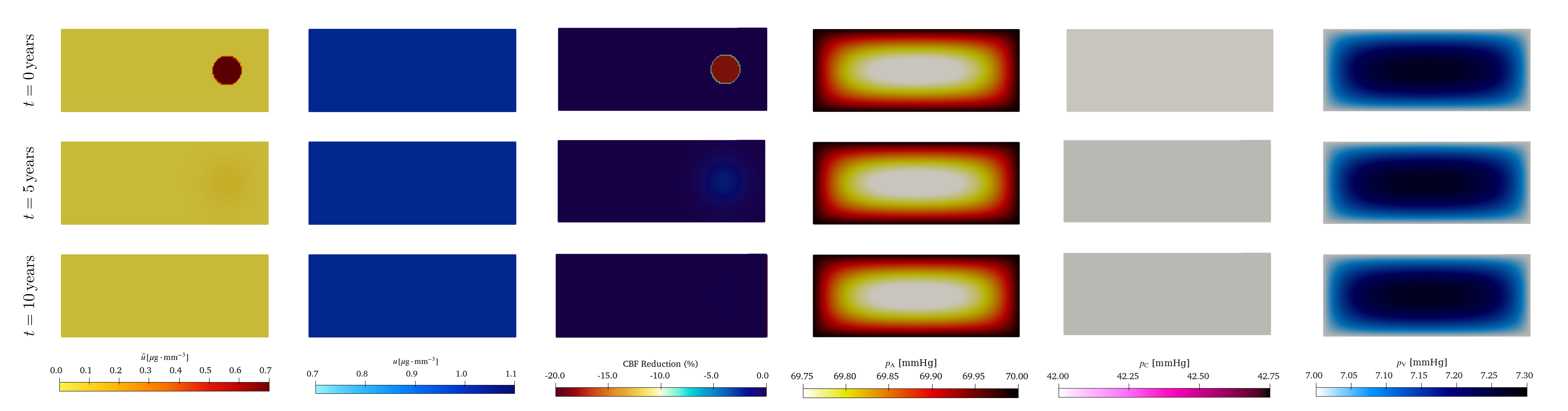}}
	\caption{Numerical solution for Test Case 1 with  a small seeding region. From left to right we have misfolded proteins $\tilde{u}$, healthy proteins $u$, reduction of CBF, arterial pressure $p_\mathrm{A}$, capillary pressure $p_\mathrm{C}$, and venous pressure $p_\mathrm{V}$. Disease propagation \textit{fails}. Values mapped to white in the color scale are shown in gray to improve background contrast.}
	\label{fig:multiple_nowavepropagation}
\end{figure}
\par
As initial data for the normal ($u_0$) and misfolded ($\tilde{u}_0$) protein concentrations, we first consider a seeding region with large radius (see Figure~\ref{fig:multiple_wavepropagation}, $t=0$):
\begin{equation*}
\text{Large:} \quad u_0(x,y)=1.0\, \mu\mathrm{g}\cdot\mathrm{mm}^{-3},\quad\tilde{u}_0(x,y)=
\begin{cases}
    0.7\,\mu\mathrm{g}\cdot\mathrm{mm}^{-3} & (x-0.08)^2+(y-0.02)^2<5\times10^{-4}, \\
    0.0\,\mu\mathrm{g}\cdot\mathrm{mm}^{-3} & \mathrm{elsewhere}.
\end{cases}
\end{equation*}
Figure~\ref{fig:multiple_wavepropagation} shows the simulation's evolution: the pathogenic protein concentration increases and saturates locally, and spreads through space like a travelling wave. Behind the wave front, the system is in a diseased state, which invades the healthy region ahead of the wave front. The CBF rate exhibits a reduction of approximately $20\%$ in the diseased region. Concerning the pressures, we observe a marked increase in capillary pressure $p_\mathrm{C}$, particularly within the pathological region but, more generally, throughout the whole domain. A similar increase is also visible in the arterial pressure $p_\mathrm{A}$, whereas the venous pressure $p_\mathrm{V}$ exhibits a decrease.

Next, we consider a seeding region with a smaller radius (see Figure~\ref{fig:multiple_nowavepropagation}, $t=0$):
\begin{equation*}
\text{Small:} \quad u_0(x,y)=1.0\,\mu\mathrm{g}\cdot\mathrm{mm}^{-3},\quad\tilde{u}_0(x,y)=
\begin{cases}
    0.7\,\mu\mathrm{g}\cdot\mathrm{mm}^{-3} & (x-0.08)^2+(y-0.02)^2<5\times10^{-5}, \\
    0.0\,\mu\mathrm{g}\cdot\mathrm{mm}^{-3} & \mathrm{elsewhere}.
\end{cases}
\end{equation*}
In contrast to the previous case, the results in Figure~\ref{fig:multiple_nowavepropagation} reveal that propagation \textit{fails}, i.e.~the initial pathogenic seed decays and the system returns to the healthy state. Moreover, no appreciable variations in the pressure distribution can be detected.

This test case demonstrates that (1) the coupled model \eqref{eq:coupled:strongformulation} has multiple stable equilibrium states, and (2) disease outbreak is dependent on the initial data. Specifically, if the initial seed of pathogenic proteins is sufficiently small (small concentration and/or small seeding region), then it decays back to the healthy equilibrium; but if the seed is large enough, then disease outbreak occurs locally and spreads like a wave through space. This behaviour is reminiscent of \textit{bistable} reaction--diffusion systems (see \cite{Keener2021}), and contrasts with most reaction--diffusion models of A\textbeta{} in the literature, which are \textit{monostable} and exhibit disease spread for every $\tilde{u}_0 \not\equiv 0$.
\subsubsection*{Numerical effects of the loosely coupled splitting strategy}
\begin{figure}[t!]
    \begin{subfigure}[b]{0.49\textwidth}
          \resizebox{\textwidth}{!}{\input{SolverComparison1.tex}}
          \caption{Comparison for the wave-propagation case.}
        \label{fig:comparison_wave}
    \end{subfigure}
    \begin{subfigure}[b]{0.49\textwidth}
        \resizebox{\textwidth}{!}{\input{SolverComparison2.tex}}
          \caption{Comparison for the wave-decay case.}
        \label{fig:comparison_decay}
    \end{subfigure}
    \caption{Test case 1: comparison between the loosely coupled splitting strategy (S, solid lines) and the fully monolithic scheme (M, dots) at the final time. The left panel shows the wave-propagation case, while the right panel reports the wave-decay case.}
    \label{fig:comparison}
\end{figure}
As discussed in Section~\ref{sec:time_discretization},  we have adopted a loosely coupled splitting strategy $(\mathbf{U}^* = \mathbf{U}^{n})$ for the solution of system~\eqref{eq:coupled:fullydiscreteformulation}. However, this choice may introduce additional dissipation and dispersion errors, which should be quantified to assess the quality of the numerical solution. For this reason, in Figure~\ref{fig:comparison} we compare the results of this test case with those obtained using a fully monolithic strategy, in which $\mathbf{U}^* = \mathbf{U}^{n+1}$ and the coupled nonlinear system is solved using a Newton method. We report the misfolded protein concentration $\tilde{u}$ along the line $\{y = 0.02\,\mathrm{m}\}$ to assess the possible presence of numerical dissipation or dispersion. In particular, Figure~\ref{fig:comparison_wave} shows the results for the case with a large initial seeding and indicates a good agreement between the waves obtained with the splitting strategy (solid lines) and those computed with the fully monolithic solver (dots). A similar behaviour is observed in the case with the small initial condition, which decays to the null equilibrium (see Figure~\ref{fig:comparison_decay}).
\subsection{Test Case 2: Injury-induced initiation of the pathology}
The following test case is motivated by the two-hit vascular hypothesis of AD (see Section~\ref{sec:intro}), which states that the initial cause of AD is vascular damage, which then triggers A\textbeta{} dyshomeostasis and disease spread \cite{Zlokovic2011,Sagare2012,kisler_cerebral_2017}. We suppose that a vascular injury affects a localized region of space, $\Omega_\mathrm{inj} \subset \Omega$, causing CBF to decrease. Focal ischaemia of this kind can arise from stroke (symptomatic or asymptomatic), small vessel disease (e.g.~in a watershed region), or atherosclerosis of an upstream artery \cite{Gorelick2011,Love2016,Snyder2015,Roher2004}.

We model the injury as a localized decrease in the transfer coefficients from the arterial to capillary compartments, representing a decrease in the supply of oxygenated arterial blood to the capillary bed. We also decrease the capillary to venous transfer coefficient and the capillary bed's permeability to model the constriction of the capillaries typically induced by those injuries \cite{Roher2004}:
\begin{equation*}
\beta_\mathrm{AC}=
\begin{cases}
    4.25\times10^{-7}\;(\mathrm{Pa\cdot s})^{-1} & \mathrm{in}\;\Omega_\mathrm{inj}, \\
    5.00\times10^{-7}\;(\mathrm{Pa\cdot s})^{-1} & \mathrm{in}\;\Omega\setminus\Omega_\mathrm{inj},
\end{cases} 
\qquad \beta_\mathrm{CV}=
\begin{cases}
    3.25\times10^{-7}\;(\mathrm{Pa\cdot s})^{-1} & \mathrm{in}\;\Omega_\mathrm{inj}, \\
    4.00\times10^{-7}\;(\mathrm{Pa\cdot s})^{-1} & \mathrm{in}\;\Omega\setminus\Omega_\mathrm{inj},
\end{cases}
\end{equation*}
\begin{equation*}
k_\mathrm{C}=
\begin{cases}
    2.00\times10^{-7}\;(\mathrm{Pa\cdot s})^{-1} & \mathrm{in}\;\Omega_\mathrm{inj}, \\
    5.00\times10^{-7}\;(\mathrm{Pa\cdot s})^{-1} & \mathrm{in}\;\Omega\setminus\Omega_\mathrm{inj}.
\end{cases}
\end{equation*}
All other parameters are spatially homogeneous and are given in Appendix A (Table \ref{tab:phys}). The initial protein concentrations are given by:
\begin{equation*}
    u_0 \equiv 1.0\,\mu\mathrm{g}\cdot\mathrm{mm}^{-3} , \quad \tilde{u}_0 \equiv 0.05\,\mu\mathrm{g}\cdot\mathrm{mm}^{-3},
\end{equation*}
representing a small perturbation from the healthy equilibrium state. The healthy CBF rate field $Q_\mathrm{H} = Q_\mathrm{H}(\boldsymbol{x})$ is the same as in Test Case 1 (Figure~\ref{fig:multiple_healthypressures}).

\begin{figure}[t!]
	\centering
	{\includegraphics[width=\textwidth]{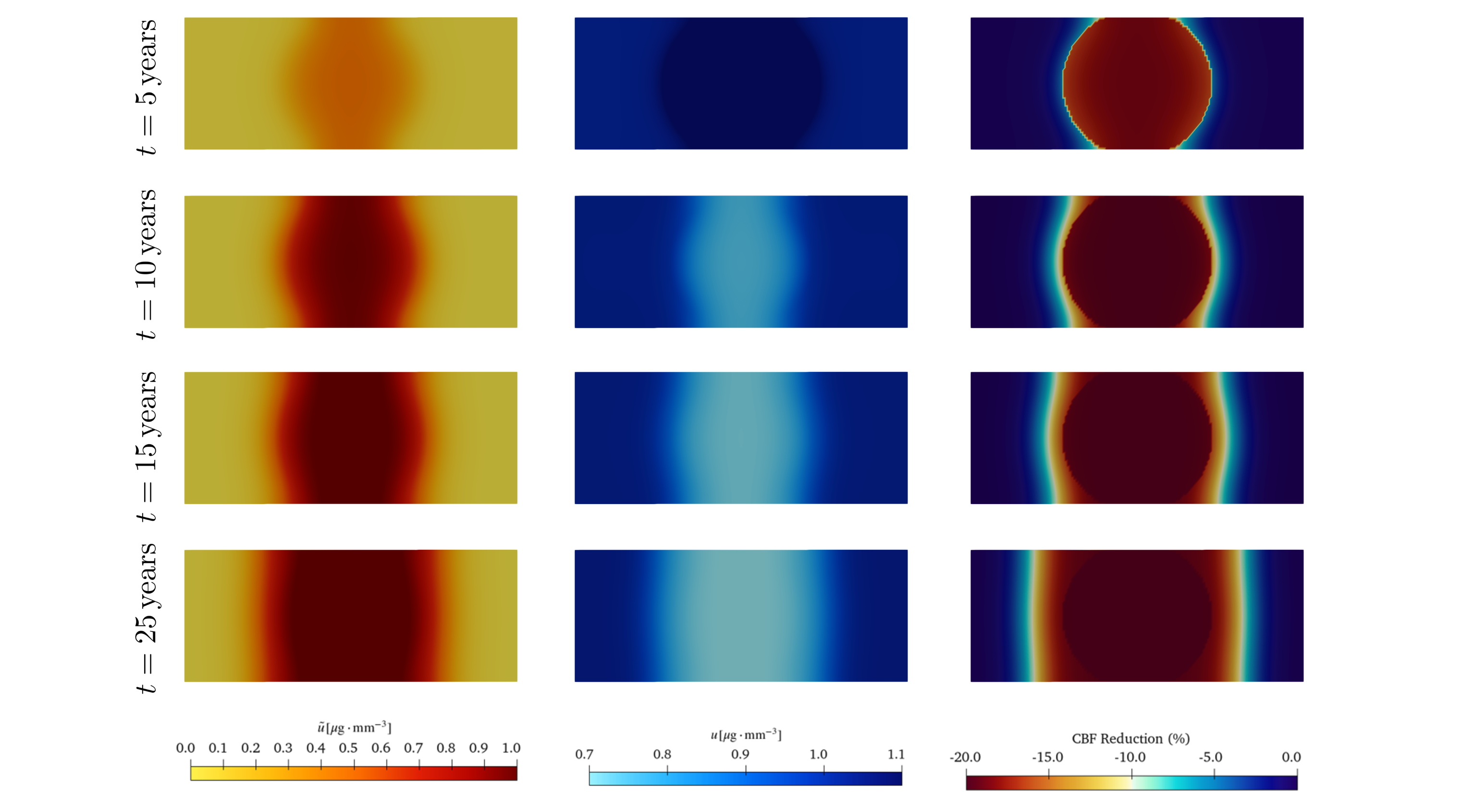}}
	\caption{Numerical solution for Test Case 2 with large injury site. From left to right we have misfolded proteins $\tilde{u}$, healthy proteins $u$, and reduction of CBF. Localized injury triggers disease spread.}
	\label{fig:hypoperfusion_wavepropagation}
\end{figure}
\begin{figure}[t!]
	\centering
	{\includegraphics[width=\textwidth]{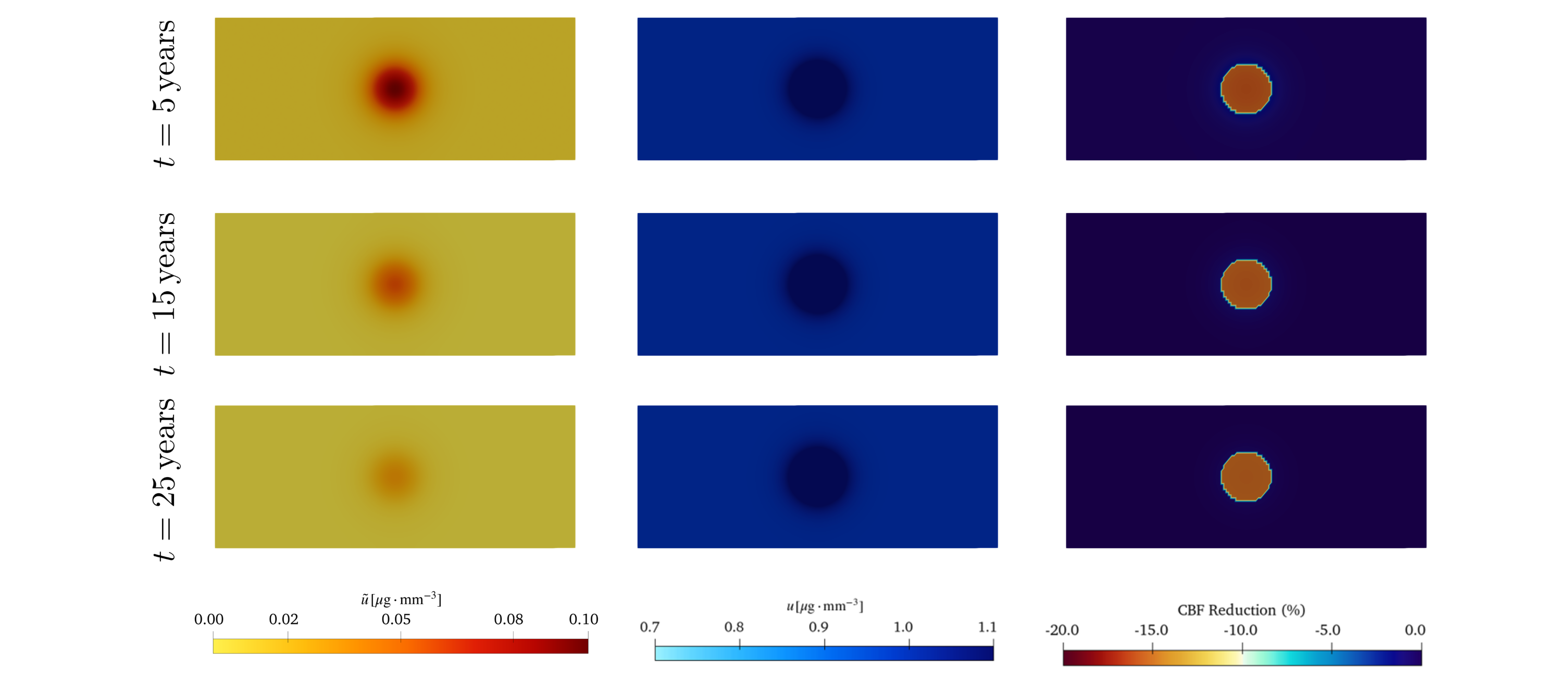}}
	\caption{Numerical solution for Test Case 2 with small injury site. From left to right we have misfolded proteins $\tilde{u}$, healthy proteins $u$, and reduction of CBF. Injury is insufficient to trigger disease spread.}
	\label{fig:hypoperfusion_nowavepropagation}
\end{figure}
\par
We first consider an injury site with large radius:
\begin{equation*}
    \text{Large injury site:} \quad \Omega_\mathrm{inj}=\{(x,y)\in\Omega \mid (x-0.05)^2+(y-0.02)^2<5\times10^{-4}\}
\end{equation*}
The simulation results are shown in Figure~\ref{fig:hypoperfusion_wavepropagation}. The reductions in $\beta_\text{AC}$, $\beta_\text{CV}$ and $k_\text{C}$ inside $\Omega_\text{inj}$ cause a $\sim$16.5\% CBF rate decrease inside $\Omega_\text{inj}$ at $t=0$. This injury-induced hypoperfusion triggers accumulation of misfolded proteins at the injury site, as seen in Figure~\ref{fig:hypoperfusion_wavepropagation}, followed by disease spread from the injury site to the rest of space. We note that the initial increase in $u$ within $\Omega_\text{inj}$ is caused by the hypoperfusion-induced increase in its production rate---see \eqref{eq:k_0^B}---and decrease in its clearance rate; after some time, $(u,\tilde{u})$ switches to its diseased state behind the invading wave front, which explains the decrease in $u$ from $t=10$ on.
\par
Next, we consider an injury site with a smaller radius:
\begin{equation*}
    \text{Small injury site:} \quad \Omega_\mathrm{inj}=\{(x,y)\in\Omega \mid (x-0.05)^2+(y-0.02)^2<5\times10^{-5}\}
\end{equation*}
In this case, we observe from the simulation results in Figure~\ref{fig:hypoperfusion_nowavepropagation} that the induced focal hypoperfusion is not sufficient to trigger disease outbreak, and the healthy state remains stable, in contrast to the case of the large injury site.

This test case is of considerable interest biologically and mathematically, as it describes the possibility of \textit{locally-induced} global disease outbreak. The intuition is that localized hypoperfusion due to vascular injury, provided it is sufficiently severe, can trigger disease outbreak locally, thus creating a bridgehead from which to invade healthy tissue; see \cite[\S5]{ahern_modelling_2025}, which introduced this idea.
\section{Numerical simulations across the whole brain}
\label{sec:num_brain}

We extend the two test cases from Section~\ref{sec:num_val} to a realistic brain geometry. The goal is to demonstrate how the properties of the mathematical model in equation~\eqref{eq:coupled:strongformulation} enable simulation of realistic pathological scenarios. Starting from a structural MRI in the OASIS-3 database~\cite{lamontagne_oasis_2019}, we generate a brain segmentation using FreeSurfer~\cite{fischl_freesurfer_2012}. We then construct a tetrahedral mesh with $142\,658$ elements (see Figure~\ref{fig:brainmeshtensors}, left) using the SVMToolkit library~\cite{mardal_mathematical_2022}.

For the simulation setup, we set a final time $T=40\,\mathrm{years}$ and time step $\Delta t=0.05\,\mathrm{years}$. Polynomial degrees are $\ell=1$ for pressures $p_\mathrm{A}$, $p_\mathrm{C}$, $p_\mathrm{V}$ and $\ell=2$ for protein concentrations $u$ and $\tilde{u}$. Physical parameters are listed in Appendix~A, Table~\ref{tab:phys_sec5}. Simulations were run on the GALILEO100 supercomputer (528 nodes, each with $2\times$Intel Cascade Lake 8260 CPUs, 24 cores, 2.4~GHz, 384~GB RAM) at the CINECA supercomputing center.

\subsubsection*{Construction of the permeability and diffusion tensors}

\begin{figure}[t]
	\centering
	{\includegraphics[width=\textwidth]{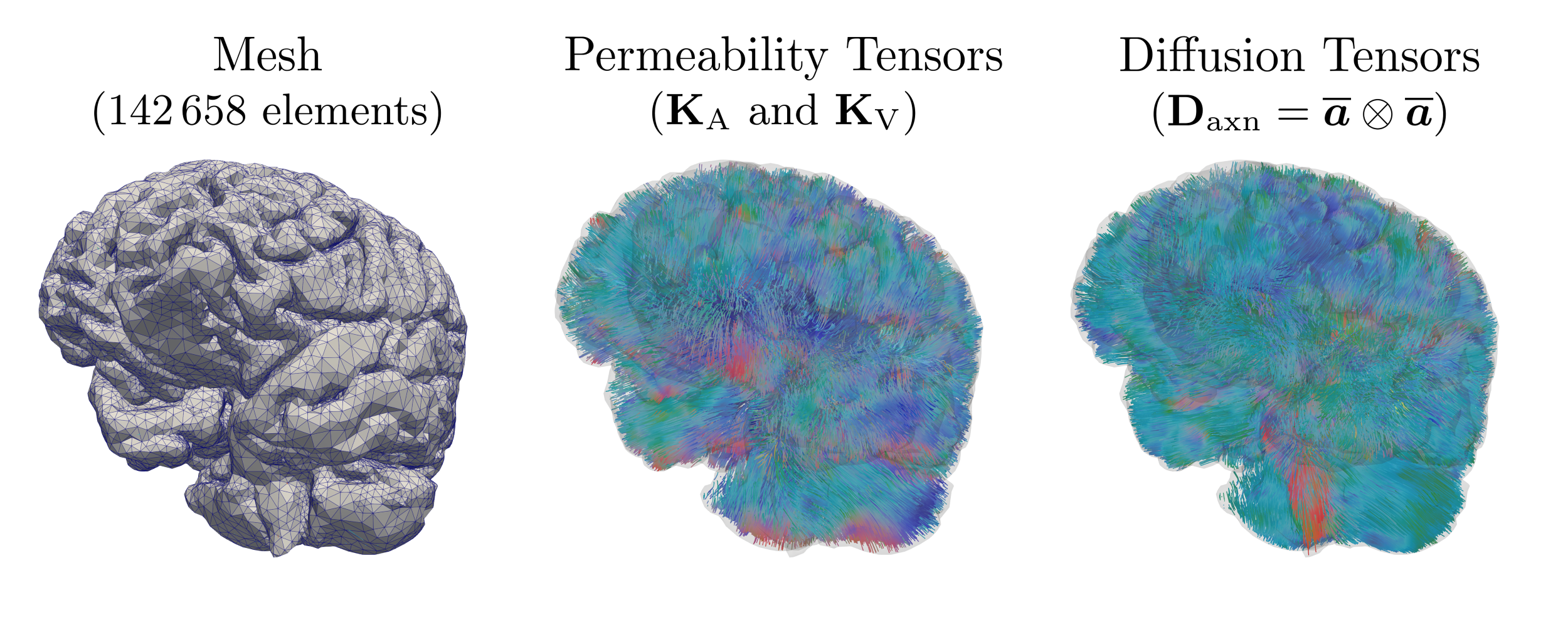}}
	\caption{Tetrahedral mesh of the brain (left), directions of permeability tensors $\mathbf{K}_\mathrm{A}$ and $\mathbf{K}_\mathrm{V}$ (centre), and axonal directions of the diffusion tensor $\mathbf{D}_\mathrm{axn}$ (right). In the visualization of the permeability and axonal directions, blue indicates directions in the $x$-axis, green indicates directions in the $y$-axis, and red indicates directions in the $z$-axis.}
	\label{fig:brainmeshtensors}
\end{figure}
To solve the perfusion porous-medium problem, we describe the perfusion tensors $\mathbf{K}_\mathrm{A}$ and $\mathbf{K}_\mathrm{V}$ for the arterioles and venules, respectively. Specifically, we follow the strategy proposed in \cite{jozsa_porous_2021} to derive directions approximately orthogonal to the pial surface, consistent with medical knowledge~\cite{schmid_depth-dependent_2017}. The resulting principal direction fibers $\overline{\boldsymbol{k}}(\boldsymbol{x})$, such that $\mathbf{K}_j = k_j (\overline{\boldsymbol{k}}\otimes \overline{\boldsymbol{k}})$ for $j=\mathrm{A,V}$, are shown in Figure~\ref{fig:brainmeshtensors} (center).

We compute the axonal component of the diffusion tensor $\mathbf{D}_\mathrm{axn}$ from Diffusion Weighted Images (DWI) using FreeSurfer~\cite{fischl_freesurfer_2012}. We then extract the principal eigenvector $\overline{\boldsymbol{a}}(\boldsymbol{x})$ to obtain the fiber directions, shown in Figure~\ref{fig:brainmeshtensors} (right). In practice, the anisotropic diffusion tensor is constructed by selecting, at each spatial location, the eigenvector associated with the largest eigenvalue of the diffusion tensor as the dominant transport direction, without taking into account the relative magnitude of the remaining eigenvalues. This diffusion model constitutes a deliberate simplification of the information provided by DWI. By prescribing a single preferred transport direction based on the principal eigenvector of the diffusion tensor, the model may introduce anisotropic transport even in regions where diffusion is nearly isotropic, and it cannot explicitly represent planar microstructural configurations characterized by comparable secondary eigenvalues. Consequently, the proposed approach should be regarded as a coarse-grained approximation of the dominant transport pathways rather than a comprehensive description of the underlying diffusion mechanisms.

\subsubsection*{Healthy CBF rate and pressures}

\begin{figure}[t]
	\centering
	{\includegraphics[width=\textwidth]{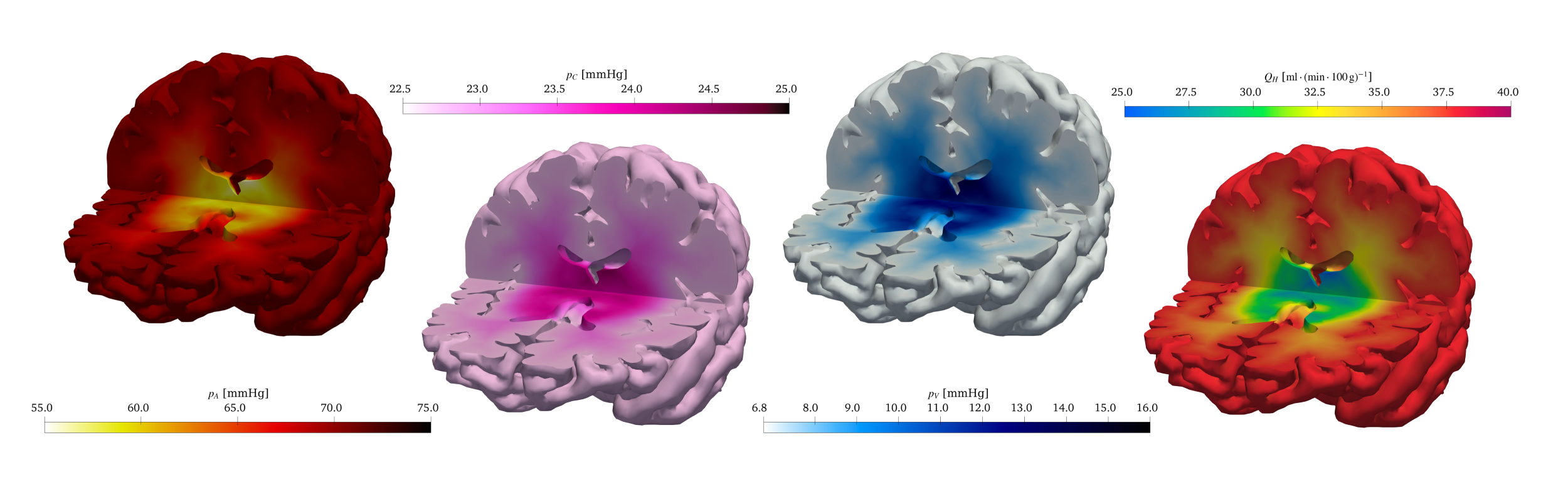}}
	\caption{Pressures computed in healthy conditions in the brain. From left to right, we report the arterial $p_\mathrm{A}$, capillary $p_\mathrm{C}$, and venous $p_\mathrm{V}$ pressures and the healthy perfusion flow $Q_\mathrm{H}$.}
	\label{fig:healthybrainpressures}
\end{figure}
As discussed in Section~4, to solve problem~\eqref{eq:coupled:strongformulation}, we first solve \eqref{eq:MP:brainproblem} for the healthy CBF rate $Q_\mathrm{H} = Q_\mathrm{H}(\boldsymbol{x})$ (in the absence of pathogenic proteins). In all brain simulations, we impose Dirichlet boundary conditions for both arterial and venous pressures on the brain's pial surface $\Gamma_\mathrm{Pial}$ and homogeneous Neumann conditions on the ventricular surface $\Gamma_\mathrm{Vent}$~\cite{corti_numerical_2023}.
\par
The solution is shown in Figure~\ref{fig:healthybrainpressures}. The arterial pressure $p_\mathrm{A}$ reaches its maximum ($\sim70\,\mathrm{mmHg}$) at the pial surface. The peaks of the capillary and venous pressures $p_\mathrm{C}$ and $p_\mathrm{V}$ occur near the center of the domain, around the ventricular wall. These computed values agree with the literature~\cite{jozsa_porous_2021}. We estimate the CBF rate inside the brain at $25$--$40\,\mathrm{ml\cdot min^{-1}\cdot (100\,\mathrm{g})^{-1}}$, consistent in magnitude with medical measurements, which show a decay in white matter relative to grey matter~\cite{fantini_cerebral_2016}.
\subsection{Test Case 3: Sensitivity to A\textbeta{-}seeding region}
\begin{figure}[t!]
    \begin{subfigure}[b]{0.49\textwidth}
          {\includegraphics[width=\textwidth]{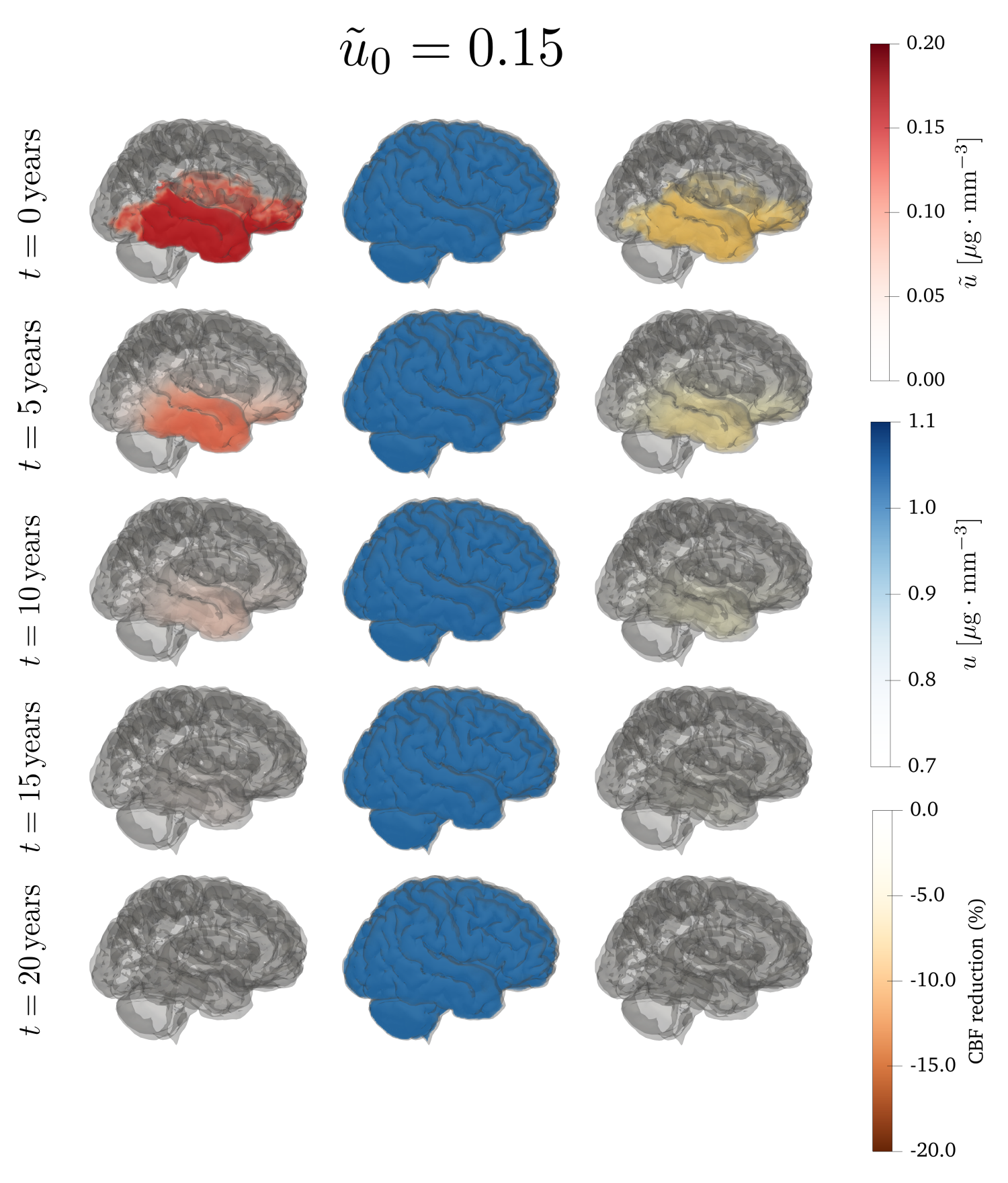}}
          \caption{Small initial seeding solution $\tilde{u}_0=0.15$.}
        \label{fig:multiple_brain_1}
    \end{subfigure}
    \begin{subfigure}[b]{0.49\textwidth}
        {\includegraphics[width=\textwidth]{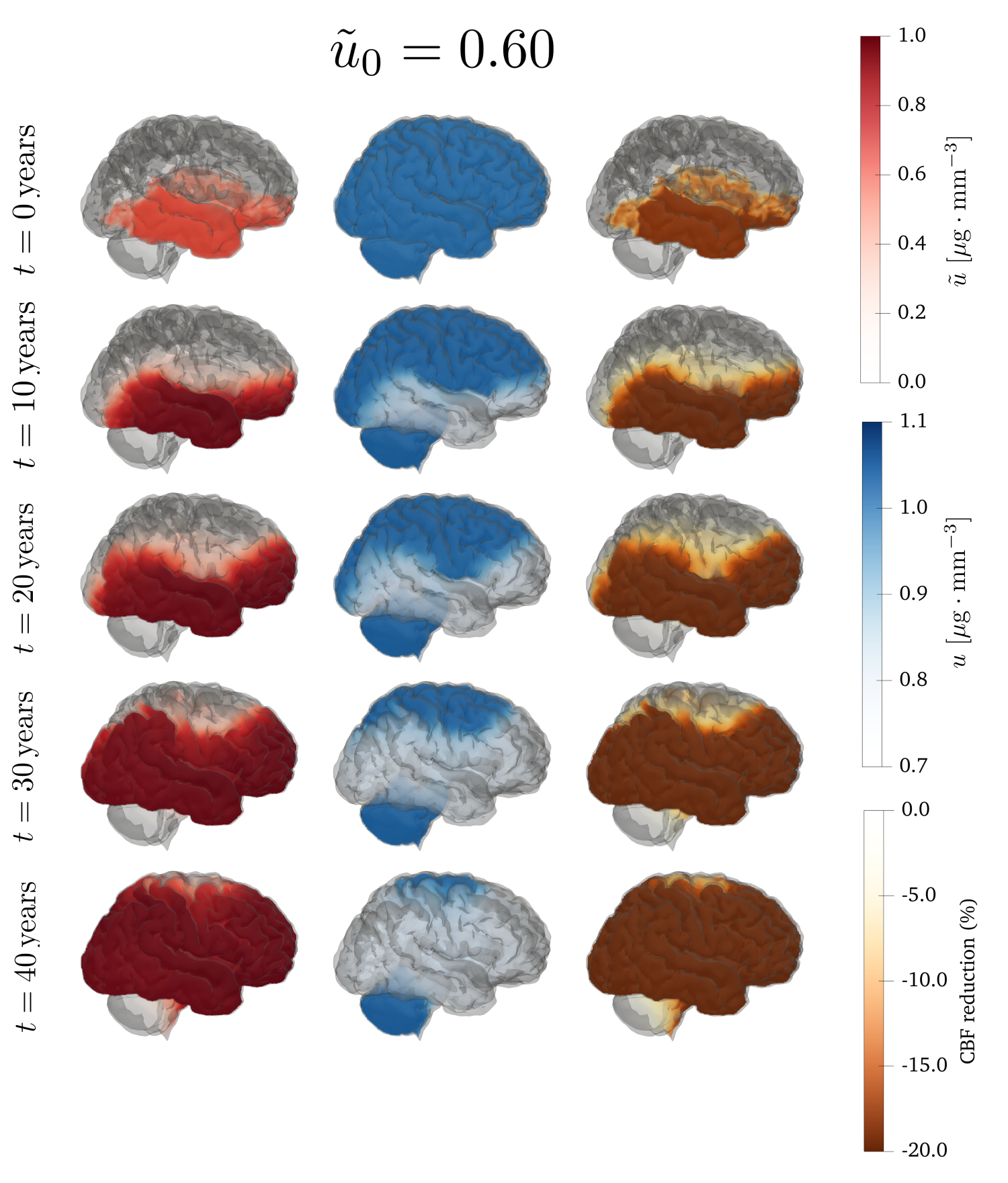}}
          \caption{Large initial seeding solution $\tilde{u}_0=0.60$.}
        \label{fig:multiple_brain_2}
    \end{subfigure}
	\caption{Test Case 3: Numerical solution comparing two different initial seeding magnitude: small seeding (a) and large seeding (b). In each panel, misfolded proteins $\tilde{u}$ (first column), healthy proteins $u$ (second column), and reduction of CBF (third column) are reported. Values mapped to white in the color scale are shown in gray to improve background contrast.}
	\label{fig:multiple_brain}
\end{figure}

As in Test Case~1, we examine the impact of the initial pathogenic A\textbeta{} distribution on the model's dynamics. Depending on its magnitude, an initial pathogenic seed can grow and spread like a wave through space or it can decay completely. Here, we simulate two initial conditions, one large and one small, both located in the basal temporal and orbitofrontal neocortex, $\Omega_\mathrm{seed}$, typical regions for the first phase of amyloid progression~\cite{goedert_alzheimer_2015}.
\par
First, we test a small localized concentration of misfolded protein ($\tilde{u}_0=0.15\,\mu\mathrm{g\cdot mm^{-3}}$ in $\Omega_\mathrm{seed}$) with constant healthy protein $u=1\,\mu\mathrm{g\cdot mm^{-3}}$ throughout the domain. Figure~\ref{fig:multiple_brain_1} shows the result. Hypoperfusion in the temporal lobe---automatically induced by our initial condition---agrees with early AD literature~\cite{eberling_reduced_1992}. The pathology does not develop: A\textbeta{} concentration decays, restoring the neocortex to the healthy equilibrium.
\par
Next, we increase the seeding magnitude to $\tilde{u}_0=0.60\,\mu\mathrm{g\cdot mm^{-3}}$ in $\Omega_\mathrm{seed}$. As shown in Figure~\ref{fig:multiple_brain_2}, this large initial seeding causes A\textbeta{} spreading through the cortical area, with associated neocortical hypoperfusion. This rise in misfolded proteins coincides with a local decrease in healthy protein population. The computed CBF rate drops by $\sim20\%$, consistent with the literature~\cite{korte_cerebral_2020}. The spreading pattern first involves the temporal lobe~\cite{marks_tau_2017}, then diffuses to the upper cortex (occipital and frontal)~\cite{thal_phases_2022}, with the parietal lobe, brainstem, and cerebellum affected last---matching the stages in~\cite{goedert_alzheimer_2015}. This reflects the spatial path length of the pathological wavefront.
\par
These simulations extend Test Case~1 to a realistic scenario with clinical relevance. They demonstrate the model's \textit{multistability}, which implies that arbitrary pathogenic seeds \textit{do not guarantee} disease accumulation and spread, in contrast to existing \textit{monostable} models of A\textbeta{} pathology \cite{weickenmeier_physics_2019,fornari_prion-like_2019}.

\subsection{Test Case 4: Injury-induced initiation of A\textbeta{} misfolding}

\begin{figure}[t!]
	\centering
	{\includegraphics[width=\textwidth]{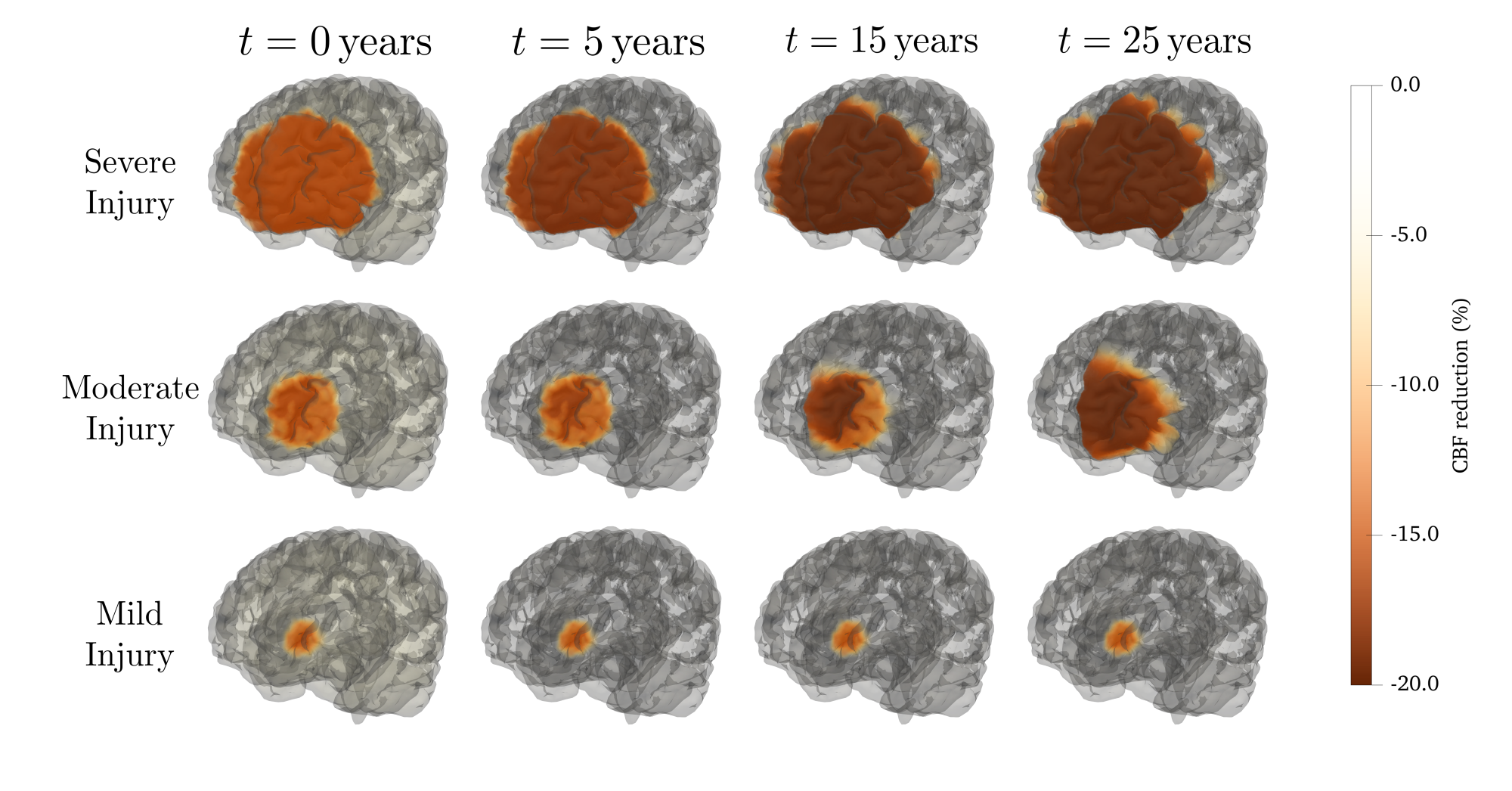}}
	\caption{Test Case 4: Reduction of CBF at different times $t=0,5,15,25$ years and for three different dimensions of the hypoperfusion region: severe (first row), moderate (second row), and mild (third row).  Values mapped to white in the color scale are shown in gray to improve background contrast.}
	\label{fig:injured_CBF}
\end{figure}
As an extension of Test Case~2, we analyze the effect of an injured subdomain size $|\Omega_\mathrm{inj}|$ on A\textbeta{} misfolding and spreading through the brain. We adopt the physical parameters from Appendix~A (Table~\ref{tab:phys_sec5}), except for
\begin{align*}
\beta_\mathrm{AC}&=
\begin{cases}
8.50\times10^{-7}\,(\mathrm{Pa\cdot s})^{-1} & \mathrm{in}\,\Omega_\mathrm{inj}, \\
1.00\times10^{-6}\,(\mathrm{Pa\cdot s})^{-1} & \mathrm{in}\,\Omega\setminus\Omega_\mathrm{inj},
\end{cases}
&
\beta_\mathrm{CV}&=
\begin{cases}
4.50\times10^{-7}\,(\mathrm{Pa\cdot s})^{-1} & \mathrm{in}\,\Omega_\mathrm{inj}, \\
3.00\times10^{-6}\,(\mathrm{Pa\cdot s})^{-1} & \mathrm{in}\,\Omega\setminus\Omega_\mathrm{inj},
\end{cases} \\[1em]
k_\mathrm{C}&=
\begin{cases}
2.50\times10^{-13}\,(\mathrm{Pa\cdot s})^{-1} & \mathrm{in}\,\Omega_\mathrm{inj}, \\
4.28\times10^{-13}\,(\mathrm{Pa\cdot s})^{-1} & \mathrm{in}\,\Omega\setminus\Omega_\mathrm{inj},
\end{cases}
\end{align*}
which induce injury in the frontal lobe of the left hemisphere. The hypoperfusion region is
$$
\Omega_\mathrm{inj}=\{(x,y,z)\in\Omega\mid(x-0.01)^2+(y-0.08)^2+(z-0.015)^2<\rho\},
$$
with radii $\rho\in\{0.01,0.02,0.04\}\,\mathrm{m}$ for mild, moderate, and severe injuries, respectively. Initial conditions are constant: $u_0=1.0\,\mathrm{\mu g}\cdot\mathrm{mm}^{-3}$ and $\tilde{u}_0=0.05\,\mathrm{\mu g}\cdot\mathrm{mm}^{-3}$ for healthy and misfolded proteins, respectively---representing a small perturbation from the healthy equilibrium. In Figure~\ref{fig:injured_CBF}, we observe that already at the initial time $t=0$, the parameters induce a significant reduction of the CBF rate in the injured region (approximately $15\%$) in all cases.
\par
\begin{figure}[t!]
	\centering
	{\includegraphics[width=\textwidth]{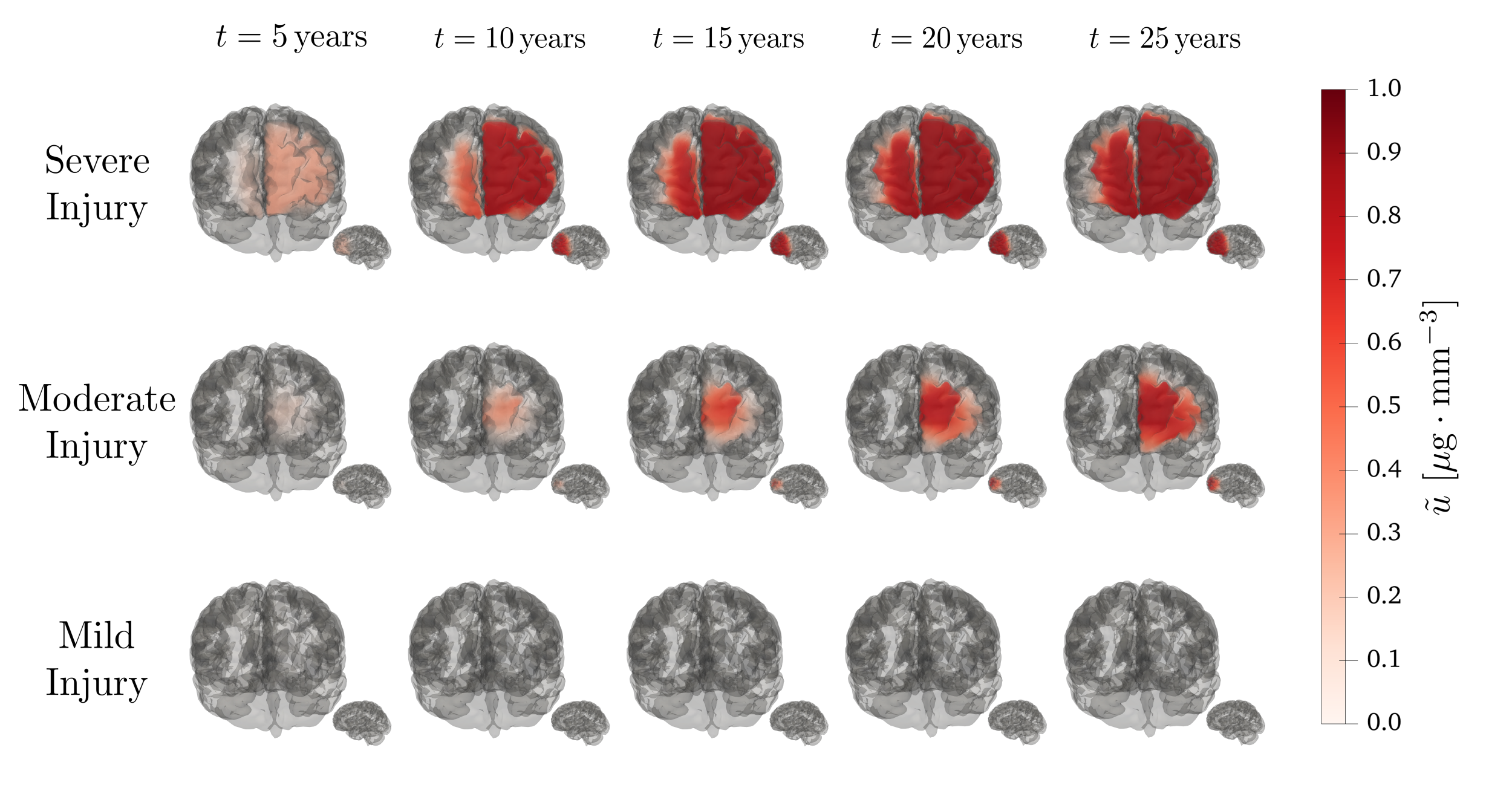}}
	\caption{Test Case 4: Misfolded protein concentration $u$ at different times $t=5,10,15,20,25$ years and for three different dimensions of the hypoperfusion region: severe (first row), moderate (second row), and mild (third row).  Values mapped to white in the color scale are shown in gray to improve background contrast.}
	\label{fig:injuried_misf}
\end{figure}
\begin{figure}[t!]
	\centering
	{\includegraphics[width=\textwidth]{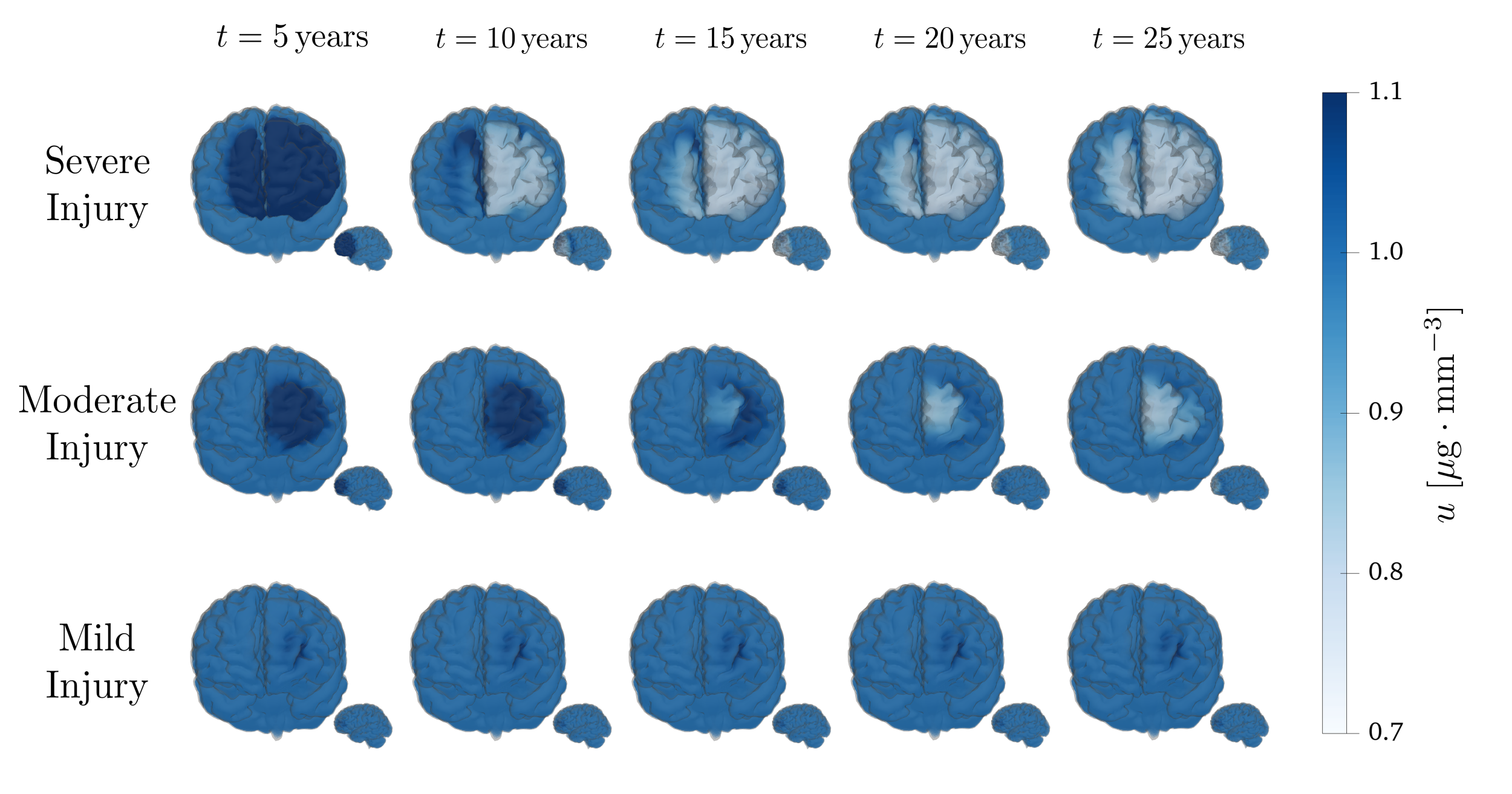}}
	\caption{Test Case 4: Healthy protein concentration $u$ at different times $t=5,10,15,20,25$ years and for three different dimensions of the hypoperfusion region: severe (first row), moderate (second row), and mild (third row).}
	\label{fig:injuried_healthy}
\end{figure}
Figures~\ref{fig:injuried_misf} and~\ref{fig:injuried_healthy} show the dynamics of misfolded and healthy protein concentrations, respectively. Moderate and severe injuries increase the pathogenic A\textbeta{} concentration ($\tilde{u}$) in the hypoperfused region, which then propagates through the frontal brain area---especially in the severe case. For the healthy A\textbeta{} concentration, we observe an initial concentration increase within the injured region, consistent with hypoperfusion-induced APP upregulation~\cite{hefter_APP_2017}, followed by a decrease in concentration as the injured region switches to the diseased equilibrium, for which $u<1$. Figure~\ref{fig:injured_CBF} shows an expansion of the hypoperfused region outward from the initial injury site in the severe and moderate cases, caused by the pathogenic A\textbeta{} distribution which induces vasoconstriction as it spreads.
\par
In the mild injury case, the pathogenic A\textbeta{} concentration decays nearly vanishes within 5 years, while healthy protein shows slight upregulation to mitigate neuronal injury~\cite{hefter_APP_2017} (see Figures~\ref{fig:injuried_misf} and~\ref{fig:injuried_healthy}). The initial hypoperfusion region remains unchanged.
\par
These results mirror those in Test Case~2. They reveal that hypoperfusive injury, if sufficiently severe, may be capable of inducing the outbreak of mixed A\textbeta{}
--vascular pathology \textit{and} its spatial invasion of otherwise healthy brain regions. This model behavior aligns with prior modelling work~\cite{ahern_modelling_2025}.
\section{Conclusion}
\label{sec:conclusion}
Alzheimer's disease, the most common cause of dementia, is characterized by the accumumulation and spread of misfolded A\textbeta{} proteins in the brain, and is also associated with chronic cerebrovascular pathology \cite{Love2016}. The work presented here is motivated by the growing recognition of a positive feedback loop between A\textbeta{} accumulation and hypoperfusion in AD \cite{Iadecola2004,Kalaria2012,korte_cerebral_2020}. Our model describes the coupling between A\textbeta{} and cerebral blood flow (CBF) associated with A\textbeta{}-induced vascoconstriction \cite{nortley_amyloid_2019,thomas_beta_1996,Niwa2001} and hypoperfusion-induced modulation of A\textbeta{} metabolism \cite{Shi1998,shi_hypoperfusion_2000,sun_hypoxia_2006,Zhang2007}.

The model is similar in spirit to that of Ahern et al.\ \cite{ahern_modelling_2025}, who assumed the same biological mechanisms. Whereas these authors used a network-based model of the brain and a semi-mechanistic model of blood flow within the network regions, the model presented here is continuous in space as well as in time, and its description of perfusion is based on multiple-network porous medium model. Interestingly, the emergent dynamics of the two models are qualitatively  similar. Both exhibit multistability (Test Cases 1 and 3) and a capacity for vascular injury-induced disease outbreak (Test Cases 2 and 4), adding to the argument that these processes emerge generically.
\par
From a numerical point of view, the differential equation has been discretized using a high-order discontinuous Galerkin method in space and implicit Euler time stepping in time. The model's multistability was demonstrated for both idealized and realistic brain geometries. Moreover, we analyzed the possibility of describing the hypoperfusion-induction of the AD by introducing an initial hypoperfusion in the mathematical model. Finally, the numerical simulations on realistic brain geometries have been performed to highlight the importance of vascular–protein coupling in the mathematical description of AD pathology. Those simulations provided the typical spreading patterns of the pathology and confirmed the model quality in describing the physical phenomena. 
\subsection{Further developments and limitations}
A limitation of the present model concerns the adoption of the heterodimer formulation, which provides a simplified description of A\textbeta{} kinetics and axonal transport. A natural extension of the present work would therefore be to enrich the reaction network by introducing additional compartments that distinguish between soluble oligomers, intermediate aggregates, and insoluble fibrils. Such multi-species models would allow us to explore state-dependent effects on clinical manifestations and comorbidities, and to represent more faithfully the different biological roles of the various aggregation states. There is no additional technical difficulty to include such effects but the real bottleneck is data availability for inference. As better data become available, it will be natural to introduce such effects.

As for the axonal transport, a purely diffusive operator cannot fully capture the directionality of axonal transport, which exhibits an anterograde–retrograde asymmetry and a concentration-dependent feedback of pathological tau.  A possible further development would be to incorporate explicit advective transport through a reformulated model with motor-dependent velocity fields, as proposed in recent network transport models~\cite{oliveri2026multiscale,tora_network-level_2025}. 

Another promising extension would be the inclusion of additional aspects of Alzheimer's disease pathology, such as brain atrophy, tau protein misfolding, and impaired clearance pathways.
Within the proposed framework, not all model parameters are individually identifiable from currently available experimental and imaging data, in particular those associated with the coupling between the A\textbeta{} and perfusion models. Nevertheless, our numerical experiments indicate that the main desired qualitative features of the dynamics are captured by the simulations. A more systematic global sensitivity analysis, aimed at quantifying the relative influence of kinetic, vascular, and coupling parameters on key outcome measures, is an important direction for future work.
\par
From a numerical perspective, an inherent drawback of the current model is its inability to automatically preserve the positivity of protein concentrations. The design of structure-preserving numerical schemes therefore represents a relevant avenue for future development. Moreover, a complete \textit{a priori} analysis of the proposed numerical method would provide further theoretical insight into its stability and convergence properties. Finally, the introduction of local discontinuous Galerkin strategies could enhance the robustness of the numerical treatment, particularly in handling the nonlinear diffusion terms.
\par
Finally, in the present formulation, anisotropy is solely based on the principal eigenvector of the diffusion tensor and does not explicitly account for tensor-shape descriptors. A more refined representation could exploit fractional anisotropy or shape indices, computed from the specific DWI data, to modulate the strength and directionality of anisotropic transport. Exploring such extensions is an interesting direction for future work.
\section*{Declaration of competing interests}
The authors declare that they have no known competing financial interests or personal relationships that could have appeared to influence the work reported in this article.

\section*{Acknowledgments}
OASIS-3 provided the brain MRI images: Longitudinal Multimodal Neuroimaging: Principal Investigators: T. Benzinger, D. Marcus, J. Morris; NIH P30 AG066444, P50 AG00561, P30 NS09857781, P01 AG026276, P01 AG003991, R01 AG043434, UL1 TR000448, R01 EB009352. AV-45 doses were provided by Avid Radiopharmaceuticals, a wholly-owned subsidiary of Eli Lilly.
\appendix
\section{Parameter values for test cases}
\label{app:param_values}
Tables~\ref{tab:phys} and \ref{tab:phys_sec5} below report the parameter values used in the simulations of Sections~\ref{sec:num_val} and \ref{sec:num_brain}, respectively.

\begin{table}[h]
	\centering
    \begin{subtable}{.5\linewidth}\centering
    {\begin{tabular}{|c|r l|}
	\hline
	   \textbf{Parameter} 
        & \multicolumn{2}{c|}{\textbf{Value}} 
        \\  \hline 
		 $\mathbf{K}_\mathrm{A}=k_\mathrm{A}\mathbf{I}$ 
        & $10^{-2}$ 
        & $[\mathrm{mm^2\cdot(Pa\cdot s)^{-1}}]$        
        \\ 
		 $\mathbf{K}_\mathrm{V}=k_\mathrm{V}\mathbf{I}$ 
        & $10^{-2}$ 
        & $[\mathrm{mm^2\cdot(Pa\cdot s)^{-1}}]$  
        \\ 
          $k_\mathrm{C}$ 
        & $5\times10^{-3}$ 
        & $[\mathrm{mm^2\cdot(Pa\cdot s)^{-1}}]$
    	\\	\hline 
          $\beta_\mathrm{AC}$ 
        & $5\times10^{-7}$ 
        & $[\mathrm{(Pa\cdot s)^{-1}}]$
        \\	
          $\beta_\mathrm{CV}$ 
        & $4\times10^{-7}$ 
        & $[\mathrm{(Pa\cdot s)^{-1}}]$
    	\\	\hline 
          $p_\mathrm{Arteries}$ 
        & $70.00$ 
        & $[\mathrm{mmHg}]$
        \\
          $p_\mathrm{Veins}$ 
        & $7.00$ 
        & $[\mathrm{mmHg}]$
    	        \\  \hline 
		 $k_\mathrm{C}^\mathrm{A\beta}$ 
        & $10^{-4}$ 
        & $[\mathrm{mm^2\cdot(Pa\cdot s)^{-1}}]$        
        \\  
	   $\alpha_{k_\mathrm{C}}$ 
        & $2.00$ 
        & $[\mathrm{years}]$ 
        \\  \hline 
		 $\beta_\mathrm{AC}^\mathrm{A\beta}$ 
        & $4\times10^{-7}$ 
        & $[\mathrm{(Pa\cdot s)^{-1}}]$        
        \\ 
	   $\alpha_{\beta_\mathrm{AC}}$ 
        & $2.00$ 
        & $[\mathrm{years}]$ 
        \\  \hline 
		 $\beta_\mathrm{CV}^\mathrm{A\beta}$ 
        & $3\times10^{-7}$ 
        & $[\mathrm{(Pa\cdot s)^{-1}}]$    
        \\  
	   $\alpha_{\beta_\mathrm{CV}}$ 
        & $2.00$ 
        & $[\mathrm{years}]$ 
    	\\	\hline 
	\end{tabular}}
    \caption{Physical parameters of perfusion equations}\label{tab:phys_perf}
    \end{subtable}%
    \begin{subtable}{.45\linewidth}\centering
    {\begin{tabular}{|c|r l|}
	\hline
	   \textbf{Parameter} 
        & \multicolumn{2}{c|}{\textbf{Value}} 
        \\  \hline 
		 $d_\mathrm{ext}$ 
        & $8.00$ 
        & $[\mathrm{mm^2\cdot years^{-1}}]$        
        \\  	 
        $d_\mathrm{axn}$ 
        & $0.00$ 
        & $[\mathrm{mm^2\cdot years^{-1}}]$        
        \\  \hline 
	   $k_1$ 
        & $1.00$ 
        & $[\mathrm{years}^{-1}]$ 
        \\  
	   $\tilde{k}_1$ 
        & $1.50$ 
        & $[\mathrm{years}^{-1}]$   
        \\  
	   $k_0$ 
        & $1.00$ 
        & $[\mu\mathrm{g}\cdot\mathrm{years}^{-1}\cdot\mathrm{mm}^{-3}]$ 
        \\  
	   $k_{12}$ 
        & $1.00$ 
        & $[\mathrm{mm}^{3}\cdot\mu\mathrm{g}^{-1}\cdot\mathrm{years}^{-2}]$ 
        \\  \hline 
		 $\kappa_1$ 
        & $1.25$ 
        & $[\mathrm{years}^{-1}]$        
        \\  	 
	   $\tilde{\kappa}_1$ 
        & $3.75$ 
        & $[\mathrm{years}^{-1}]$   
        \\  
	   $\kappa_0$ 
        & $1.25$ 
        & $[\mu\mathrm{g}\cdot\mathrm{years}^{-1}\cdot\mathrm{mm}^{-3}]$ 
        \\  \hline 
	\end{tabular}}
    \caption{Physical parameters of heterodimer equations}\label{tab:phys_heterodimer}
    \end{subtable} 
        \\[6pt]
    \centering
    \begin{subtable}{.5\linewidth}
    \centering
    \begin{tabular}{|c|c|c|}
	\hline
	   \textbf{Dimensionless} 
        & \textbf{Healthy} & \textbf{Misfolded} \\	   \textbf{Parameter} 
        & \textbf{Value} & \textbf{Value} 
       \\  \hline
		 $\boldsymbol{\sigma}_\mathrm{A}=\sigma_\mathrm{A}\left(\overline{\boldsymbol{k}}\otimes \overline{\boldsymbol{k}}\right)$ 
        & \multicolumn{2}{c|}{$3.00\times10^{4}$}      
        \\ 
		 $\boldsymbol{\sigma}_\mathrm{V}=\sigma_\mathrm{V}\left(\overline{\boldsymbol{k}}\otimes \overline{\boldsymbol{k}}\right)$ 
        & \multicolumn{2}{c|}{$3.00\times10^{4}$} 
        \\
        $B$ 
        & \multicolumn{2}{c|}{$0.75$} 
    	\\ \hline
          $\boldsymbol{\sigma}_\mathrm{C}=\sigma_\mathrm{C}\mathbf{I}$ 
        & $1.50\times10^{4}$ 
        & $3.00\times10^{2}$ 
    	\\ 
          $\tilde{\gamma}_{\mathrm{CV}}$ 
        & $1.00$ & $0.75$ 
    	\\ 
          $\tilde{\gamma}_{\mathrm{AC}}$ 
        & $1.00$ & $0.80$
    	\\	\hline 
	\end{tabular}
    \caption{Dimensionless parameters of perfusion equations}\label{tab:phys_dimless_perf}
    \end{subtable}%
        \begin{subtable}{.45\linewidth}\centering
    {\begin{tabular}{|c|c|c|}
	\hline
	   \textbf{Dimensionless} 
        & \textbf{Healthy} & \textbf{Misfolded} \\	   \textbf{Parameter} 
        & \textbf{Value} & \textbf{Value}
       \\  \hline
		 $\delta_\mathrm{ext}$
        & \multicolumn{2}{c|}{$1.00$}         
        \\ 
		 $\delta_\mathrm{axn}$
        & \multicolumn{2}{c|}{$0.00$}         
        \\
        $\epsilon$ 
        & \multicolumn{2}{c|}{$1.50$}
        \\ 
          $R$ 
        & \multicolumn{2}{c|}{$0.66$} 
    	\\ \hline
        $\lambda^\mathrm{B}_1$ 
        & $1.00$ & $0.75$
    	\\
        $\tilde{\lambda}^\mathrm{B}_1$ 
        & $1.00$ & $0.50$
    	\\
        $ \mu^\mathrm{B}_0$ 
        & $1.00$ & $0.75$
    	\\
        \hline 
	\end{tabular}}
    \caption{Dimensionless parameters of heterodimer equations}\label{tab:phys_dimless_het}
    \end{subtable}%
	\caption{Test Cases of Section 4: Physical parameters appearing in Equation~\eqref{eq:coupled:strongformulation}. The values are arbitrarily chosen for the simulation on the idealized geometry.}
	\label{tab:phys}
\end{table}
\begin{table}[ht!]
	\centering
    \begin{subtable}{.5\linewidth}\centering
    \begin{tabular}{|c|r l|}
	\hline
	   \textbf{Parameter} 
        & \multicolumn{2}{c|}{\textbf{Value}} 
       \\  \hline
		 $\mathbf{K}_\mathrm{A}=k_\mathrm{A}\left(\overline{\boldsymbol{k}}\otimes \overline{\boldsymbol{k}}\right)$ 
        & $1.23\times10^{-3}$ 
        & $[\mathrm{mm^2\cdot(Pa\cdot s)^{-1}}]$        
        \\ 
		 $\mathbf{K}_\mathrm{V}=k_\mathrm{V}\left(\overline{\boldsymbol{k}}\otimes \overline{\boldsymbol{k}}\right)$ 
        & $1.23\times10^{-3}$ 
        & $[\mathrm{mm^2\cdot(Pa\cdot s)^{-1}}]$  
        \\ 
          $\mathbf{K}_\mathrm{C}=k_\mathrm{C}\mathbf{I}$ 
        & $4.28\times10^{-7}$ 
        & $[\mathrm{mm^2\cdot(Pa\cdot s)^{-1}}]$
    	\\	\hline 
          $\beta_\mathrm{AC}$ 
        & $1.00\times10^{-6}$ 
        & $[\mathrm{(Pa\cdot s)^{-1}}]$
        \\	
          $\beta_\mathrm{CV}$ 
        & $3.00\times10^{-6}$ 
        & $[\mathrm{(Pa\cdot s)^{-1}}]$
    	\\	\hline 
          $p_\mathrm{Arteries}$ 
        & $70.00$ 
        & $[\mathrm{mmHg}]$
        \\
          $p_\mathrm{Veins}$ 
        & $7.00$ 
        & $[\mathrm{mmHg}]$
    	\\  \hline 
		 $k_\mathrm{C}^\mathrm{A\beta}$ 
        & $1.00\times10^{-7}$ 
        & $[\mathrm{mm^2\cdot(Pa\cdot s)^{-1}}]$        
        \\  
	   $\alpha_{k_\mathrm{C}}$ 
        & $2.50$ 
        & $[\mathrm{years}]$ 
        \\  \hline 
		 $\beta_\mathrm{AC}^\mathrm{A\beta}$ 
        & $8.00\times10^{-7}$ 
        & $[\mathrm{(Pa\cdot s)^{-1}}]$        
        \\ 
	   $\alpha_{\beta_\mathrm{AC}}$ 
        & $2.50$ 
        & $[\mathrm{years}]$ 
        \\  \hline 
		 $\beta_\mathrm{CV}^\mathrm{A\beta}$ 
        & $2.40\times10^{-6}$ 
        & $[\mathrm{(Pa\cdot s)^{-1}}]$    
        \\  
	   $\alpha_{\beta_\mathrm{CV}}$ 
        & $2.50$ 
        & $[\mathrm{years}]$ 
    	\\	\hline 
	\end{tabular}
    \caption{Physical parameters of perfusion equations.}\label{tab:phys_perf_sec5}
    \end{subtable}%
    \begin{subtable}{.45\linewidth}\centering
    \begin{tabular}{|c|r l|}
	\hline
	   \textbf{Parameter} 
        & \multicolumn{2}{c|}{\textbf{Value}} 
        \\  \hline 
		 $d_\mathrm{ext}$ 
        & $8.00$ 
        & $[\mathrm{mm^2\cdot years^{-1}}]$        
        \\  	 
        $d_\mathrm{axn}$ 
        & $80.00$ 
        & $[\mathrm{mm^2\cdot years^{-1}}]$        
        \\  \hline 
	   $k_1$ 
        & $1.00$ 
        & $[\mathrm{years}^{-1}]$ 
        \\  
	   $\tilde{k}_1$ 
        & $1.50$ 
        & $[\mathrm{years}^{-1}]$   
        \\  
	   $k_0$ 
        & $1.00$ 
        & $[\mu\mathrm{g}\cdot\mathrm{years}^{-1}\cdot\mathrm{mm}^{-3}]$ 
        \\  
	   $k_{12}$ 
        & $1.20$ 
        & $[\mathrm{mm}^{3}\cdot\mu\mathrm{g}^{-1}\cdot\mathrm{years}^{-2}]$ 
        \\ \hline 
		 $\kappa_1$ 
        & $1.25$ 
        & $[\mathrm{years}^{-1}]$        
        \\  	 
	   $\tilde{\kappa}_1$ 
        & $3.75$ 
        & $[\mathrm{years}^{-1}]$   
        \\  
	   $\kappa_0$ 
        & $1.25$ 
        & $[\mu\mathrm{g}\cdot\mathrm{years}^{-1}\cdot\mathrm{mm}^{-3}]$ 
        \\  \hline 
	\end{tabular}
    \caption{Physical parameters of heterodimer equations.}\label{tab:phys_blood_heterodimer_sec5}
    \end{subtable}%
    \\[6pt]
    \centering
    \begin{subtable}{.5\linewidth}
    \centering
    \begin{tabular}{|c|c|c|}
	\hline
	   \textbf{Dimensionless} 
        & \textbf{Healthy} & \textbf{Misfolded} \\	   \textbf{Parameter} 
        & \textbf{Value} & \textbf{Value} 
       \\  \hline
		 $\boldsymbol{\sigma}_\mathrm{A}=\sigma_\mathrm{A}\left(\overline{\boldsymbol{k}}\otimes \overline{\boldsymbol{k}}\right)$ 
        & \multicolumn{2}{c|}{$2.31\times10^{2}$}      
        \\ 
		 $\boldsymbol{\sigma}_\mathrm{V}=\sigma_\mathrm{V}\left(\overline{\boldsymbol{k}}\otimes \overline{\boldsymbol{k}}\right)$ 
        & \multicolumn{2}{c|}{$2.31\times10^{2}$} 
        \\
        $B$ 
        & \multicolumn{2}{c|}{$3.00$} 
    	\\ \hline
          $\boldsymbol{\sigma}_\mathrm{C}=\sigma_\mathrm{C}\mathbf{I}$ 
        & $8.02\times10^{-2}$ 
        & $1.87\times10^{-2}$ 
    	\\ 
          $\tilde{\gamma}_{\mathrm{CV}}$ 
        & $1.00$ & $0.80$ 
    	\\ 
          $\tilde{\gamma}_{\mathrm{AC}}$ 
        & $1.00$ & $0.80$
    	\\	\hline 
	\end{tabular}
    \caption{Dimensionless parameters of perfusion equations.}\label{tab:phys_dimless_perf_sec5}
    \end{subtable}%
        \begin{subtable}{.45\linewidth}\centering
    {\begin{tabular}{|c|c|c|}
	\hline
	   \textbf{Dimensionless} 
        & \textbf{Healthy} & \textbf{Misfolded} \\	   \textbf{Parameter} 
        & \textbf{Value} & \textbf{Value}
       \\  \hline
		 $\delta_\mathrm{ext}$
        & \multicolumn{2}{c|}{$1.00$}         
        \\ 
		 $\delta_\mathrm{axn}$
        & \multicolumn{2}{c|}{$10.00$}         
        \\
        $\epsilon$ 
        & \multicolumn{2}{c|}{$1.50$}
        \\ 
          $R$ 
        & \multicolumn{2}{c|}{$0.80$} 
    	\\ \hline
        $\lambda^\mathrm{B}_1$ 
        & $1.00$ & $0.75$
    	\\
        $\tilde{\lambda}^\mathrm{B}_1$ 
        & $1.00$ & $0.50$
    	\\
        $ \mu^\mathrm{B}_0$ 
        & $1.00$ & $0.75$
    	\\
        \hline
	\end{tabular}}
    \caption{Dimensionless parameters of heterodimer equations.}\label{tab:phys_dimless_het_sec5}
    \end{subtable}%
	\caption{Test Cases of Section 5: Physical and dimensionless parameters of the model. The values of the parameters of the uncoupled problems of (a) from \cite{jozsa_porous_2021} and of (b) from \cite{antonietti_discontinuous_2024}.}
    
	\label{tab:phys_sec5}
\end{table}

\bibliographystyle{abbrv}
\bibliography{bibliography}

\newpage

\end{document}